\documentclass{article}

\usepackage[english]{babel}
\usepackage[letterpaper,top=2cm,bottom=2cm,left=3cm,right=3cm,marginparwidth=1.75cm]{geometry}

\usepackage{amsmath, bm, amsthm, bbm}
\usepackage{graphicx}
\usepackage[colorlinks=true, allcolors=blue]{hyperref}
\usepackage{subfig}

\usepackage{float}
\usepackage{amsfonts}

\usepackage{xcolor}

\usepackage{algorithm}
\usepackage{algpseudocode}

\usepackage{multirow} 

\newtheorem{theo}{Theorem}[section]
\newtheorem{remark}[theo]{Remark}
\newtheorem{ex}[theo]{Example}

\usepackage{soul}

\usepackage{lineno}

\title{A general class of continuous asymmetric distributions with positive support
}
\usepackage{authblk} 
\author[1]{Felipe S. Quintino  \thanks{Corresponding author: felipes.quintino2@gmail.com}}
\author[1]{Pushpa N. Rathie  \thanks{pushpanrathie@yahoo.com}}
\author[2]{Luan C. S. M. Ozelim \thanks{luanoz@gmail.com}}
\author[3]{Tiago A. da Fonseca \thanks{fonsecafga@unb.br}}
\author[1]{Roberto Vila \thanks{rovig161@gmail.com}}

 \affil[1]{\normalfont 
	Department of Statistics, University of
	Bras\'ilia, Bras\'ilia 70.910-900, Brazil}
  \affil[2]{\normalfont 
	Department of Civil and Environmental Engineering, University of
	Bras\'ilia, Bras\'ilia 70.910-900, Brazil}
 \affil[3]{\normalfont 
	Gama Engineering College, University of
	Bras\'ilia, Bras\'ilia 72.444-240,  Brazil}

\date{\today}

\begin{document}
\maketitle
\begin{abstract}
In order to better fit real-world datasets, studying asymmetric distribution is of great interest. 
In this work, we derive several mathematical properties of a general class of asymmetric distributions with positive support which shows up as a unified framework for Extreme Value Theory asymptotic results. 
The new model generalizes some well-known distribution models such as Generalized Gamma, Inverse Gamma, Weibull, Fréchet, Half-normal, Modified half-normal, Rayleigh, and Erlang.
To highlight the applicability of our results, the performance of the analytical models is evaluated through real-life dataset modeling.

\textbf{Keywords}:  
{ Generalized classes of distributions, Skewed probability distributions, Extreme value $\mathbb{H}$-function, Maximum likelihood estimator, Statistical modeling of asymmetric data.}
\end{abstract}

\section{Introduction}

Asymmetric distributions are common across various fields, such as epidemiology and finance, where they help model phenomena like disease spread and risk assessment \cite{huang2019epirank}. In finance, heavy-tailed distributions are the most suitable tool for predicting stock performance and optimizing portfolios by quantifying probabilities and managing risk \cite{Rathie2017,Rathie2016,quintino2024estimation}. Our focus, however, is to study a new class of distributions with positive support that generalizes many existing asymmetric models.

Nadarajah \cite{Nadarajah03} studied problems involving stress-strength reliability (SSR) for extreme value distributions. He wrote the probability $P(X<Y)$ in terms of generalized hypergeometric functions when $X$ and $Y$ were independent random variables with Weibull, Fréchet, or Gumbel distribution and relied on severe parameter restrictions to get analytical results.

Inspired by the work of \cite{Nadarajah03}, \cite{Rathieetal2023} 
observed that the parameter restrictions could be relaxed if a new class of special functions was introduced, called extreme value $\mathbb{H}$-function. Later, \cite{oliveira2024stress,quintino2024asset,quintino2024estimation} showed that $\mathbb{H}$-functions could also be used to describe SSR for the distributions of Generalized Extreme Value (GEV) distributions, transmuted GEV and p-max stable laws families.

Departing from mathematical properties of $\mathbb{H}$-functions studied in \cite{Rathieetal2023}, a new class of probability density function (PDF) with positive support was introduced 
{ in \cite{Rathieetal2023}. However, the mathematical and statistical properties and the applications of this model require a more in-depth study, which will be provided in this work. The PDF is given by:}
\begin{equation}\label{eq_intr_g}
	g(y; \boldsymbol{\theta}) = \frac{1}{c(\boldsymbol{\theta})} y^{\theta_6}\exp\left(-\theta_1y-(\theta_2y^{\theta_3}+\theta_4)^{\theta_5}\right), ~~y>0,
\end{equation}  
where $\boldsymbol{\theta}=(\theta_1,\theta_2, \theta_3, \theta_4, \theta_5, \theta_6)$ and $c(\boldsymbol{\theta})$ is a normalizing constant (more details about the parameter space for $\boldsymbol{\theta}$ will be given in the following sections).
This new class of PDFs generalizes several well-known models in the literature, such as exponential, Gamma, Generalized Gamma \cite{stacy1962generalization}, Inverse Gaussian \cite{johnson1995continuous}, Inverse Generalized Gamma \cite{ramos2021bayesian}, Weibull \cite{johnson1995continuous}, Fréchet, Half-normal, Modified half-normal \cite{sun2023modified}, Rayleigh \cite{johnson1995continuous}, Maxwell–Boltzmann \cite{johnson1995continuous} and Erlang \cite{johnson1995continuous} distributions, among others.

{ Although \cite{Rathieetal2023} introduced the PDF \eqref{eq_intr_g}, they focused only on the mathematical properties of the extreme value $\mathbb{H}-$function, thus a complete characterization of the model needs to be discussed to fully expose its potential. For example, a topic that needs to be studied is the description of the incomplete extreme value $\mathbb{H}-$function, which is necessary for presenting the model's cumulative distribution function (CDF). The stochastic representation provides an alternative way to generate random variables, which was not addressed in the previous work. A more in-depth study of the modality and shapes, which are essential for a proper understanding of the model’s characteristics, is also needed to be described, in addition to other important features such as moments and characteristic function. These aspects are fundamental not only for ensuring the internal consistency of the proposed probability model but also for enabling its practical use in statistical inference and simulation-based analysis. Without a complete characterization, including the CDF, stochastic representation, and moment properties, the applicability of the distribution to real-world problems remains limited. Therefore, this work aims to bridge this gap by providing both the theoretical development and the computational tools required for its implementation.}

The aim of this paper is to study statistical properties of \eqref{eq_intr_g}, such as particular cases, modality, shapes, finite mixtures, characteristic function, Mellin transform, moments, and entropy. Inferential properties are also addressed through maximum likelihood estimators and conditions for $g(y;\boldsymbol{\theta})$ to belong to the exponential family are derived. A second estimator is also proposed, based on the difference between the empirical distribution function and the theoretical cumulative distribution function.

In our study, we assess two estimators through real data modeling. We examine three distinct datasets to gain insights into their performance.
The first dataset focuses on the minimum monthly water flows (measured in cubic meters per second, m³/s) of the Piracicaba River, situated in São Paulo state, Brazil.
The second dataset investigates the strength of carbon fibers under stress. These fibers were tested to understand their behavior when subjected to tension.
Finally, we analyze a dataset that captures the failure times of various machine parts. This information is important for maintenance and reliability considerations.

The paper is organized as follows: In Section \ref{sec_preliminaries} the model is presented and we study modality, shapes, and particular cases.
In Section \ref{sec_main}, properties of the new class of distributions are derived. Two classes of estimators are presented in Section \ref{sec_estimation}, while Section \ref{sec_applications} presents real data applications.

\section{A class of continuous distributions with positive support}\label{sec_preliminaries}

Let $\Theta\subset\mathbb{R}^6$ be a parameter space.
For a random variable $Y$ we define the PDF 
\begin{equation}\label{eq_pdf}
  g(y; \boldsymbol{\theta}) = \frac{1}{c(\boldsymbol{\theta})}  y^{\theta_6}\exp\left( -\theta_1 y -(\theta_2 y^{\theta_3} + \theta_4)^{\theta_5} \right), ~~y>0,
\end{equation}
where $\boldsymbol{\theta}=(\theta_1,\theta_2,\theta_3,\theta_4,\theta_5,\theta_6)\in \Theta$ is a parameter vector such that $\theta_1\geq0$, $\theta_2\geq0$ ($\theta_1$ and $\theta_2$ cannot be zero simultaneously), $\theta_4\geq0$ and $\theta_3, \theta_5\in\mathbb{R}$. 
$\theta_6>-1$ when $\theta_1\neq0$ or $\theta_1=0$ and ${\rm sign}(\theta_3)={\rm sign}(\theta_5)$, $\theta_6<-1$ when $\theta_1=0$ and ${\rm sign}(\theta_3)\neq {\rm sign}(\theta_5)$. These are general restrictions which could be further manipulated when $\theta_5 \in \mathbb{Z}$. In particular, an integer valued $\theta_5$ allows one to relax the constraints on $\theta_1$, $\theta_2$ and $\theta_4$, extending their possible domains to negative values. This should, on the other hand, be assessed in a case by case basis. In Subsection \ref{ap_partialH}, for example, the case when $\theta_5 \in \mathbb{N}$ is discussed. In that case, if $\theta_5$ is even, then $\theta_1>0$ and  $\theta_2,\theta_4 \in \mathbb{R}$. On the other hand, if $\theta_5$ is odd, then three possible scenarios arise: if $\theta_3 \theta_5<1$, then $\theta_2>0$  and $\theta_1 \in \mathbb{R}$; if  $\theta_3 \theta_5=1$, then $(\theta_1 + \theta_2^{\theta_5})>0$; finally, if $\theta_3 \theta_5 <1$, then $\theta_1>0$ and $\theta_2 \in \mathbb{R}$.
The function $c({\boldsymbol{\theta}})$ appearing in the definition of $g$ is a normalizing factor written as
\begin{equation}\label{eq_H_function}
    c(\boldsymbol{\theta}) = \mathbb{H}(\boldsymbol{\theta}) := \mathbb{H}(\theta_1,\theta_2,\theta_3,\theta_4,\theta_5,\theta_6) = \int_0^\infty y^{\theta_6}\exp\left( -\theta_1 y -(\theta_2 y^{\theta_3} + \theta_4)^{\theta_5} \right) {\rm d}y.
\end{equation}
We recall that extreme value $\mathbb{H}$-function was first defined in \cite{Rathieetal2023} who also proposed the PDF \eqref{eq_pdf}. Statistical and mathematical properties of $g(y; \boldsymbol{\theta})$ are addressed in this paper.

{ Note that, in general, the model \eqref{eq_pdf} is not (globally) identifiable. For example, for $\boldsymbol{\theta}=(\theta, 0,0,0,1,0)\neq \boldsymbol{\theta}'=(0,\theta,1,0,1,0)$, we have $g(y;\boldsymbol{\theta}) = g(y;\boldsymbol{\theta}')$. By globally identifiable, we mean that the mapping $\Theta \ni \boldsymbol{\theta} \mapsto g(y; \boldsymbol{\theta})$ is one-to-one, for all $y > 0$, onto the entire parameter space $\Theta$. It is often possible to achieve identifiability in a parametric model that fails to satisfy this property by introducing specific technical constraints, such as those discussed in Ref. \cite{rothenberg1971identification}, thereby facilitating meaningful inference.}

\begin{figure}[H]
    \centering
    \includegraphics[width=1.0\linewidth]{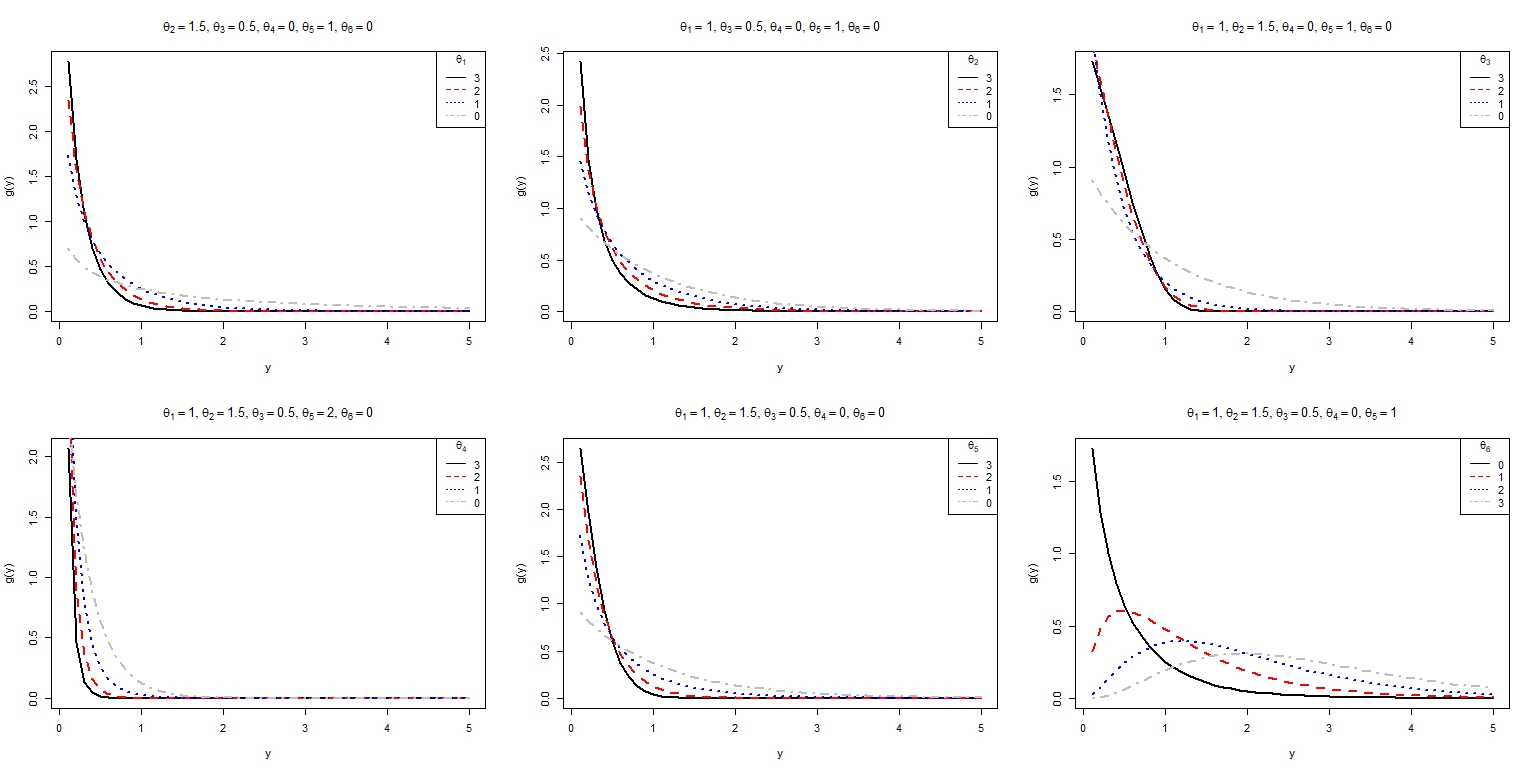}
    \caption{Plots of the PDF $g(y; \boldsymbol{\theta})$.}
    \label{fig:pdf}
\end{figure}

The CDF corresponding to PDF $g(y; \boldsymbol{\theta})$ is given by
\begin{equation}\label{eq_CDF}
    G(x) = G(x; \bm \theta) := \frac{1}{c(\boldsymbol{\theta})}\mathbb{H}(x;\boldsymbol{\theta}) = \frac{\mathbb{H}(x;\boldsymbol{\theta})}{\mathbb{H}(\boldsymbol{\theta})}, ~~x>0,
\end{equation}
where $\mathbb{H}(x;\boldsymbol{\theta})=\int_0^x y^{\theta_6}\exp\left( -\theta_1 y -(\theta_2 y^\theta_3 + \theta_4)^{\theta_5} \right) {\rm d}y $ is the incomplete extreme value $\mathbb{H}$-function (for more details see Subsection \ref{ap_partialH}). Figure \ref{fig:pdf} presents the behavior the PDFs of this random variable for some arbitrary values of its parameters.

\begin{remark}
Note that the PDF \eqref{eq_pdf} of $g(y; \boldsymbol{\theta}) $ defines a weighted gamma distribution by considering the gamma distribution with shape parameter $\theta_6+1$ and rate parameter $\theta_1$ and the weight function $w(y) = \exp(-(\theta_2 y^{\theta_3} + \theta_4)^{\theta_5})$, $y>0$.
\end{remark}

\subsection{The incomplete extreme value { $\mathbb{H}$-function}}\label{ap_partialH}

Consider the extreme value $\mathbb{H}$-function in \eqref{eq_H_function}. Define 
$$h_y(\boldsymbol{\theta}) := y^{\theta_6} \exp\left(-\theta_1 y - (\theta_2 y^{\theta_3}+\theta_4)^{\theta_5}\right), ~y>0.$$
Then, the incomplete extreme value $\mathbb{H}-$function is defined by
\begin{equation}\label{ap:H_partial}
    \mathbb{H}(x;\boldsymbol{\theta}) := \int_0^x y^{\theta_6}\exp\left( -\theta_1 y -(\theta_2 y^{\theta_3} + \theta_4)^{\theta_5} \right) {\rm d}y, ~x>0.
\end{equation}

Observe that $h_y(\boldsymbol{\theta}) = c(\boldsymbol{\theta})  g(y; \boldsymbol{\theta})$. Then, the CDF corresponding to $g(y; \boldsymbol{\theta})$ is given by
$$G(x;\boldsymbol{\theta}) = \frac{1}{c(\boldsymbol{\theta})}\int_0^x h_y(\boldsymbol{\theta}) {\rm d}y, ~x>0. $$

When $\theta_5=m\in\mathbb{N}$, it follows from Eq.(16) in \cite{Rathieetal2023} that
\begin{equation}\label{ap:H_function_m}
    \mathbb{H}(\theta_1, \theta_2, \theta_3, \theta_4, m, \theta_6) = \sum_{n=0}^\infty \frac{(-1)^n \theta_4^{mn}}{n!} \sum_{k=0}^{mn}\left(\begin{array}{c}
         mn \\
         k
    \end{array}\right) \left(\frac{\theta_2}{\theta_4}\right)^k \frac{\Gamma(\theta_6+\theta_3 k+1)}{\theta_1^{\theta_6+\theta_3 k+1}}.
\end{equation}

\noindent thus, we can write \eqref{ap:H_partial} in terms of series as follows:
\begin{equation}\label{ap:H_incomplete}
      \mathbb{H}( x;\theta_1, \theta_2, \theta_3, \theta_4, m, \theta_6) = \sum_{n=0}^\infty \frac{(-1)^n \theta_4^{mn}}{n!} \sum_{k=0}^{mn}\left(\begin{array}{c}
           mn  \\
           k
      \end{array}\right) \left(\frac{\theta_2}{\theta_4}\right)^k \frac{\gamma(\theta_6+\theta_3k+1,x)}{\theta_1^{\theta_6+\theta_3 k+1}}, ~x>0,
\end{equation}
where $\gamma(p,x)$ is the incomplete gamma function $\gamma(p,x)=\int_{0}^{x}w^{p-1}{\rm e}^{-w} {\rm d}w.$
Note that by taking $x\to\infty$ in \eqref{ap:H_incomplete} we get \eqref{ap:H_function_m}.

Another important function is $\int_x^\infty h_y(\boldsymbol{\theta}){\rm d}y$. It can be used, for example, to define the survival function and the Hazard function. Thus, 
\begin{eqnarray}
    \nonumber \int_x^\infty h_y(\boldsymbol{\theta}){\rm d}y &=& \mathbb{H}(\theta_1, \theta_2, \theta_3, \theta_4, m,  \theta_6) - \mathbb{H}( x;\theta_1, \theta_2, \theta_3, \theta_4, m, \theta_6),
\end{eqnarray}
which can be further manipulated, by considering \eqref{ap:H_function_m} and \eqref{ap:H_incomplete}, to obtain
\begin{eqnarray}
    \nonumber \int_x^\infty h_y(\boldsymbol{\theta}){\rm d}y &=& 
    \sum_{n=0}^\infty \frac{(-1)^n \theta_4^{mn}}{n!} \sum_{k=0}^{mn}\left(\begin{array}{c}
           mn  \\
           k
      \end{array}\right) \left(\frac{\theta_2}{\theta_4}\right)^k \left[ \frac{\Gamma(\theta_6+\theta_3 k+1)-\gamma(\theta_6+\theta_3k+1,x)}{\theta_1^{\theta_6+\theta_3 k+1}}\right], ~~x>0.
\end{eqnarray}

\subsection{Particular probability models}

\begin{figure}[H]
    \centering
        \includegraphics[width=1.0\linewidth]{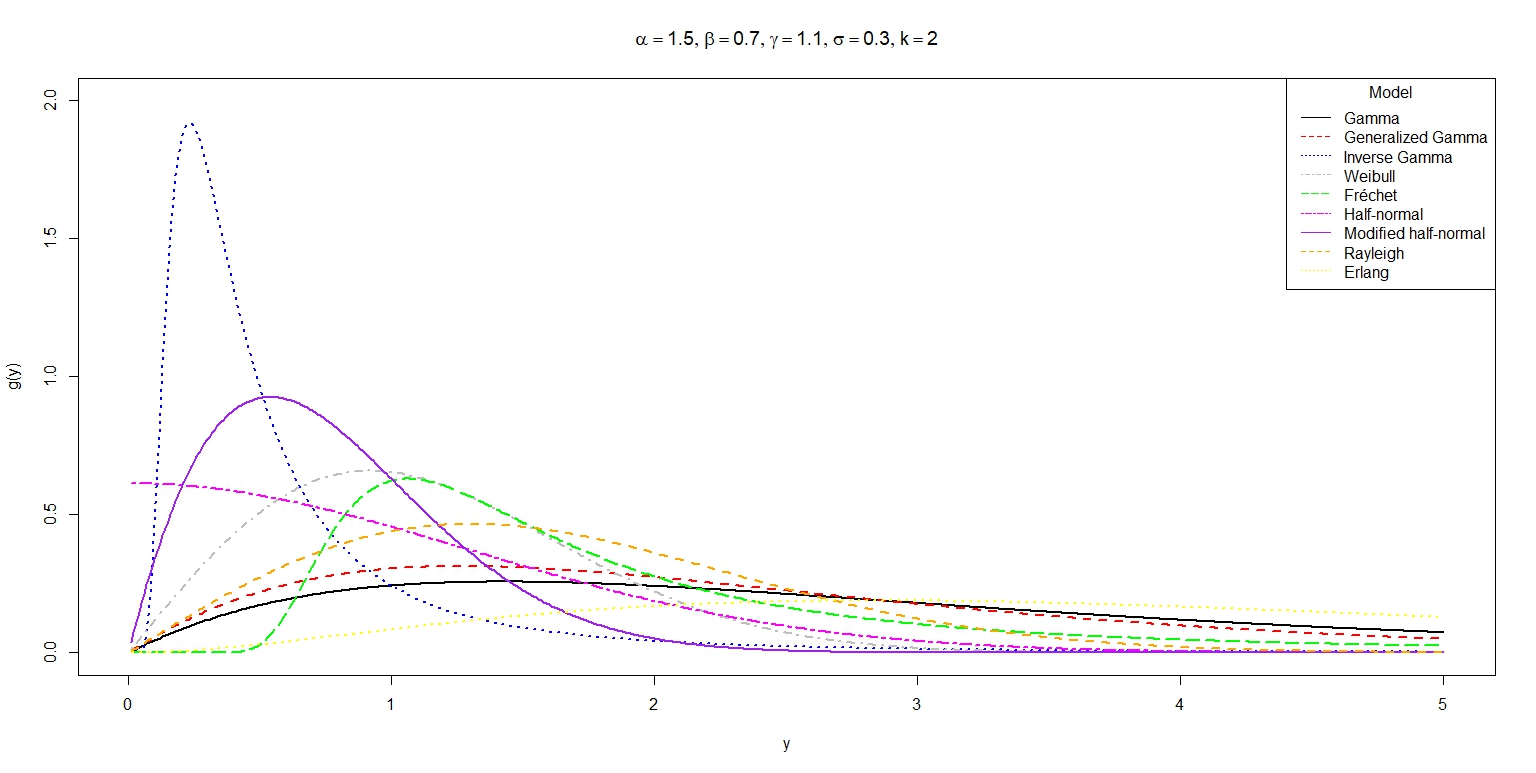}

    \caption{Plots of the PDF $g(y; \boldsymbol{\theta})$ for particular models presented in Table \ref{tab:g_particular}.}
    \label{fig:pdf_particular}
\end{figure}

Some particular examples of densities $g(y; \boldsymbol{\theta})$ are presented in Table \ref{tab:g_particular}. Figure \ref{fig:pdf_particular} presents the plots of the PDFs of such particular cases. There are other models that are not particular cases of \eqref{eq_pdf}, but it can be written as a finite mixture of $g(y; \boldsymbol{\theta})$ (see Subsection \ref{subsec_mixture}).

\begin{table}[H]
\caption{Particular cases of PDF $g(y; \boldsymbol{\theta})$.}
\label{tab:g_particular}
\centering
\begin{tabular}{lcc}
 \hline
Distribution     & PDF & $\boldsymbol{\theta}$   \\
\hline
Gamma { \cite{johnson1995continuous}}& $\frac{\beta^\alpha}{\Gamma(\alpha)}y^{\alpha-1}\exp\left(-\beta y\right)$, $\alpha, \beta>0$   & $(\beta,0,1,0,1,\alpha-1)$ \\
Generalized Gamma { \cite{stacy1962generalization}} & $\frac{\gamma\beta^{\alpha/\gamma}}{\Gamma(\alpha/\gamma)}y^{\alpha-1}\exp\left(-\beta y^\gamma\right)$, $\alpha, \beta,\gamma>0$    & $(0,\beta,\gamma,0,1,\alpha-1)$ \\
Inverse Gamma { \cite{ramos2021bayesian}}& $\frac{\beta^{\alpha}}{\Gamma(\alpha)}y^{-(\alpha+1)}\exp\left(-\beta \frac{1}{y}\right)$, $\alpha, \beta>0$      & $(0,\beta,-1,0,1,-\alpha-1)$ \\
Weibull { \cite{johnson1995continuous}}&  $\frac{\alpha}{\sigma}\left(\frac{y}{\sigma}\right)^{\alpha-1}\exp\left(-\left(\frac{y}{\sigma} \right)^\alpha\right)$, $\alpha, \sigma>0$    & $(0,1/\sigma,1,0,\alpha,\alpha-1)$  \\
Fréchet { \cite{Ramos2020Frechet}}&  $\frac{\alpha}{\sigma}\left(\frac{y}{\sigma}\right)^{-\alpha-1}\exp\left(-\left(\frac{y}{\sigma} \right)^{-\alpha}\right)$, $\alpha, \sigma>0$    & $(0,1/\sigma,1,0,-\alpha,-\alpha-1)$  \\
Half-normal { \cite{johnson1995continuous}}&   $\frac{\sqrt{2}}{\sigma \sqrt{\pi}} \exp\left( -\frac{y^2}{2\sigma^2}\right)$,  $ \sigma>0$    & $(0,\frac{1}{2\sigma^2},2,0,1,0)$ \\
Modified half-normal { \cite{sun2023modified}}&    $\frac{2\beta^{\alpha/2}}{\Psi\left(\alpha/2, \gamma/\sqrt{\beta} \right)}y^{\alpha-1}\exp\left(\gamma x - \beta x^2 \right)$, $\alpha, \beta>0,\gamma<0$ & $(-\gamma,\beta,2,0,1,\alpha-1)$ \\
Rayleigh { \cite{johnson1995continuous}}&  $\frac{1}{\sigma^2}y\exp\left(-\frac{1}{2\sigma^2}y^2 \right)$, $\sigma>0$    & $(0,\frac{1}{2\sigma^2},2,0,1,1)$ \\
Erlang { \cite{johnson1995continuous}}&   $\frac{\beta^k}{(k-1)!}y^{k-1}\exp\left(-\beta y \right)$, $\beta>0$, $k=1,2,\cdots$   & $(\beta,0,1,1,1,k-1)$  \\
\hline
\end{tabular}
\end{table}

It is interesting to highlight that Gumbel random variables are a limiting case of the new class of random variables hereby presented. In that case, one may notice that:

\begin{equation}
    \exp(z) = \lim_{n\to\infty}\left(1+\frac{z}{n}\right)^n 
\end{equation}

Thus, considering the PDF of Gumbel distribution with location parameter $\alpha \in \mathbb{R}$ and scale $\beta>0$:

\begin{eqnarray}
      e^{-y/\beta+ \alpha/\beta- e^{-y/\beta +\alpha/\beta}}/\beta &=& e^{-y/\beta+ \alpha/\beta- \lim_{n\to\infty}\left(1+\frac{-y/\beta +\alpha/\beta}{n}\right)^n}/\beta \nonumber \\
     &=& \lim_{n\to\infty} e^{-y/\beta+ \alpha/\beta- \left(1+ \frac{\alpha}{n\beta} - \frac{y}{\beta n}\right)^n}/\beta \\
     &=& e^{\alpha/\beta} \lim_{n\to\infty} h_y(1/ \beta,-1/ \beta n,1,1+\alpha/ \beta n,n,0)/\beta\nonumber
\end{eqnarray}

\subsection{Stochastic representation}

Let $Y$ be a random variable with PDF $g(\cdot; \boldsymbol{\theta})$ given in \eqref{eq_pdf}.
In this subsection, we prove that $Y$ has the following stochastic representation
\begin{align}\label{rep-st}
    Y\stackrel{d}{=}\sum_{n=0}^\infty
    X_n 1_{\{W_n=n\}},
\end{align}
where ``$\stackrel{d}{=}$'' denotes equality in distribution, $X_n$ and $W_n$ are independent, $X_n$ is an absolutely continuous random variable with CDF (for $x>0$) 
\begin{align}\label{id-1}
    F_{X_n}(x)
    =
    \dfrac{E[(\theta_2 Z^{\theta_3}+\theta_4)^{n\theta_5} {1}_{\{Z\leq x\}}]}{E[(\theta_2 Z^{\theta_3}+\theta_4)^{n\theta_5}]},\quad  Z\sim{\rm Gamma(\theta_6+1,\theta_1)}
\end{align}
and
$W_n$ is a discrete random variable with probability function
\begin{align}\label{id-2}
P(W_n=n)={(-1)^n\over n!}\, \dfrac{\theta_6\Gamma(\theta_6) E[(\theta_2 Z^{\theta_3}+\theta_4)^{n\theta_5}]}{c(\boldsymbol{\theta})\theta_1^{\theta_6+1}},
    \quad n=0,1,\ldots.
\end{align}

    Indeed, by using the law of total probability and by \eqref{rep-st} (for $y>0$),
    \begin{align}\label{id-3}
        P(Y\leq y)
        &=
        \sum_{k=0}^\infty
        P(Y\leq y\vert W_k=k) P(W_k=k) \nonumber
        \\[0,2cm]
        &=
        \sum_{k=0}^\infty
        P(X_k\leq y\vert W_k=k) P(W_k=k) \nonumber
               \\[0,2cm]
        &=
        \sum_{k=0}^\infty
        F_{X_k}(y) P(W_k=k),
    \end{align}
    where in the last line the independence of $X_n$ and $W_n$  was used. By employing \eqref{id-1} and \eqref{id-2} in \eqref{id-3} the CDF of $Y$ is written as (for $y>0$)
    \begin{align*}
        P(Y\leq y)
        &=
        \dfrac{1}{c(\boldsymbol{\theta})}
        \sum_{k=0}^\infty
        {(-1)^k\over k!}\, 
        E[(\theta_2 Z^{\theta_3}+\theta_4)^{k\theta_5} {1}_{\{Z\leq y\}}],
        \quad Z\sim{\rm Gamma(\theta_6+1,\theta_1)},
                    \\[0,2cm]
        &=
        \dfrac{\theta_6\Gamma(\theta_6)}{c(\boldsymbol{\theta}) \theta_1^{\theta_6+1}}\,
        E\left[\exp\{-(\theta_2 Z^{\theta_3}+\theta_4)^{\theta_5}\} {1}_{\{Z\leq y\}}\right]
                            \\[0,2cm]
        &=
                \dfrac{1}{c(\boldsymbol{\theta})}\,
\int_0^y
s^{\theta_6}\exp\left( -\theta_1 s -(\theta_2 s^{\theta_3} + \theta_4)^{\theta_5} \right) {\rm d}s
=
\dfrac{\mathbb{H}(y;\boldsymbol{\theta}) 
}{\mathbb{H}(\boldsymbol{\theta})}
=
G(y),
    \end{align*}
    where in the last identity we have used \eqref{eq_CDF}.
    We have thus completed the proof of \eqref{rep-st}.

    \begin{remark}
        It is clear that when $\theta_5\in\mathbb{N}$ the expectations in \eqref{id-1} and \eqref{id-2} become simple closed-forms depending on the CDF of the gamma distribution.
    \end{remark}

\subsection{Modality and shapes}

A simple calculation shows that
\begin{align*}
\frac{d}{dy}
g(y; \boldsymbol{\theta}) 
= 
\frac{g(y; \boldsymbol{\theta})}{c(\boldsymbol{\theta}) y} 
\left[
\theta_6
-
\theta_1 y
-
\theta_2\theta_3\theta_5
(\theta_2 y^{\theta_3} + \theta_4)^{\theta_5-1}
y^{\theta_3}
\right]=0
\end{align*}
iff
\begin{align}\label{eq-1}
\theta_6
-
\theta_1 y
=
\theta_2\theta_3\theta_5
(\theta_2 y^{\theta_3} + \theta_4)^{\theta_5-1}
y^{\theta_3}.
\end{align}
Therefore, every critical point $y$ of $g$ is a positive solution of the above equation.

In what follows we analyze the number of roots of Equation \eqref{eq-1} in the less trivial cases.
\begin{enumerate}
    \item 
If [$\theta_2=0$ or $\theta_3=0$ or $\theta_5=0$], $\theta_1>0$ and $\theta_6>0$ then $y=\theta_6/\theta_1$ is the unique solution of \eqref{eq-1}.
    \item 
If $-1<\theta_6\leq 0$, [$\theta_1\neq 0$ or $\theta_1= 0$] and ${\rm sign}(\theta_3)={\rm sign}(\theta_5)$  then it is clear that \eqref{eq-1} has no solution.
    \item 
If $\theta_6<-1$, $\theta_1=0$, ${\rm sign}(\theta_3)\neq{\rm sign}(\theta_5)$ and $\theta_2, \theta_4>0$ then \eqref{eq-1} it is equivalent to
\begin{align}\label{eq-2}
\frac{\theta_6}{\theta_2\theta_3\theta_5
}
=
(\theta_2 z + \theta_4)^{\theta_5-1}
z, 
\quad z=y^{\theta_3}.
\end{align}
First, let us assume that $\theta_5>0$. In this case,
the function $z\mapsto r(z)=(\theta_2 z + \theta_4)^{\theta_5-1}z$ is increasing, and hence,  \eqref{eq-2} has a unique positive solution.
On the other hand, assume that $\theta_5<0$.
In this case, a routine calculation shows that $z\mapsto r(z)$ is unimodal, with maximum point 
$z_0=-\theta_4/(\theta_2\theta_5)$ and maximum value 
\begin{align*}
r(z_0)
=
\left(\frac{1}{-\theta_5} + 1\right)^{\theta_5-1}
\frac{\theta_4^{\theta_5}}{\theta_2(-\theta_5)}.
\end{align*}
Imposing the condition
$
\left({1/(-\theta_5)} + 1\right)^{\theta_5-1}
\theta_4^{\theta_5}
>
{-\theta_6/ \theta_3
}
$
we have $r(z_0)>\theta_6/(\theta_2\theta_3\theta_5)$, and therefore \eqref{eq-2} has two positive solutions. Which implies that \eqref{eq-1} also has two positive solutions.
On the other hand, under condition 
$
\left({1/(-\theta_5)} + 1\right)^{\theta_5-1}
\theta_4^{\theta_5}
<
{-\theta_6/ \theta_3
}
$
we have $r(z_0)<\theta_6/(\theta_2\theta_3\theta_5)$,  and therefore \eqref{eq-2} has no positive solutions. Which implies that \eqref{eq-1} also has no positive solutions.
    \item 

If $\theta_6>0$, [$\theta_1\neq 0$ or $\theta_1= 0$] and ${\rm sign}(\theta_3)={\rm sign}(\theta_5)$  then the function $y\mapsto s(y)=(\theta_2 y^{\theta_3} + \theta_4)^{\theta_5-1}
y^{\theta_3}$ is increasing when $\theta_5>0$ and is unimodal when $\theta_5<0$. It is clear that in the case $\theta_5>0$ the equation \eqref{eq-1} has a unique positive root. Following the steps in the previous item, it can be proven that, when $\theta_5<0$ the equation \eqref{eq-1} has two or no positive solutions.
\end{enumerate}

Since $\lim_{y\to\infty}g(y;\boldsymbol{\theta})=0$, the number of roots deduced in Items 1-4 above show that the density $g$ in \eqref{eq_pdf} has strictly decreasing, unimodal and decreasing-increasing-decreasing forms.

It is also possible to build a general solution to (\ref{eq-1}) where one of the roots can be expressed by using Lagrange's Inversion Theorem \cite{Lagrange}. In that case, let $y$ be defined as the following function of constant $\chi$, function $\phi$, and a parameter $\delta$:

\begin{equation}\label{eq-LIT1}
    y = \chi + \delta \phi(y)
\end{equation}

\noindent then any function $\zeta(y)$ is expressed as the following power series in $\delta$:

\begin{equation}\label{eq-LIT2}
    \zeta(y) = \zeta(\chi) + \sum_{n=1}^{\infty} \frac{\delta^n}{n!} \frac{d^{n-1}}{dx^{n-1}} \left\{ \frac{d \zeta(x)}{dx}\phi^n(x) \right\} \bigg \rvert_{x = \chi},
\end{equation}

\noindent where suitable restrictions apply to make sure the series above converges. This way, equation (\ref{eq-1}) can be rearranged as:

\begin{equation}\label{eq-LIT3}
y
=
\frac{\theta_6}{\theta_1 } - 
\frac{\theta_2\theta_3\theta_5}{\theta_1 }
(\theta_2 y^{\theta_3} + \theta_4)^{\theta_5-1}
y^{\theta_3}.
\end{equation}

The direct comparison of (\ref{eq-LIT1}), (\ref{eq-LIT2}) and (\ref{eq-LIT3}) reveals that, for $\zeta(y)=y$:

\begin{eqnarray}\label{eq-LIT4}
    y &=& \frac{\theta_6}{\theta_1 } + \sum_{n=1}^{\infty} \frac{(-1)^n \theta_2^n\theta_3^n\theta_5^n}{\theta_1^n n!} \frac{d^{n-1}}{dx^{n-1}} \left\{(\theta_2 x^{\theta_3} + \theta_4)^{\theta_5 n-n}
x^{\theta_3 n} \right\} \bigg \rvert_{x = \frac{\theta_6}{\theta_1 }} \nonumber \\
 &=& \frac{\theta_6}{\theta_1 } - \frac{\theta_2\theta_3\theta_5 \theta_6^{\theta_3}}{\theta_1^{\theta_3+1} } \left(\frac{\theta_2 \theta_6^{\theta_3}}{\theta_1^{\theta_3} } + \theta_4 \right)^{\theta_5-1}+ \cdots 
\end{eqnarray}

\noindent under suitable constraints.

\subsection{Finite mixtures}\label{subsec_mixture}
Some recently introduced distributions can be written as a finite mixture of $g$ densities. These are the cases of a bimodal Weibull distribution--introduced by \cite{vila2022bimodal}-- and the transmuted GEV distribution --introduced by \cite{aryal2009transmuted}.

It follows from \eqref{eq_pdf} and Eq.(6) in \cite{vila2022bimodal} that the PDF of the bimodal Weibull can be written as:
\begin{eqnarray}
   \nonumber f(y;\alpha,\beta,\delta) &=& c_0(\alpha,\beta,\delta)g(y; (0,1/\beta, 1,0,\alpha,\alpha-1)) + c_1(\alpha,\beta,\delta)g(y; (0,1/\beta, 1,0,\alpha,\alpha))\\
   \nonumber &+& c_2(\alpha,\beta,\delta)g(y; (0,1/\beta, 1,0,\alpha,\alpha+1)),
\end{eqnarray}
where $c_j=cj(\alpha,\beta,\delta)$ are constants depending only on $(\alpha, \beta, \delta)$ and $c_0+c_1+c_2=1$.

The PDF of TGEV model depends on $(\mu, \sigma, \gamma, \lambda)\in\mathbb{R}\times\mathbb{R}_+\times\mathbb{R}\times[-1,1]$. The support of the model is $(\mu-\sigma/\gamma,\infty)$, if $\gamma>0$, or $(-\infty, \mu-\sigma/\gamma)$, if $\gamma<0$, or $\mathbb{R}$, if $\gamma=0.$
We are interested in the case where TGEV is used to model positive data. In this case, we consider $\gamma>0$ and $\mu=\sigma/\gamma$. It follows from \eqref{eq_pdf} and Eq.(4) in \cite{aryal2009transmuted} that the PDF of TGEV distribution is 
\begin{eqnarray}
    \nonumber f_2(y;\sigma, \gamma, \lambda) &=& c_0 g(y; (0,\gamma/\sigma, 1, 0,-1/\gamma, -1-1/\gamma)) + c_1 g(y; (0,2^{-\gamma} \frac{\gamma}{\sigma}, 1, 0,-1/\gamma, -1-1/\gamma)), 
\end{eqnarray}
where $c_0$ and $c_1$ are constants depending only on $(\sigma, \gamma, \lambda)$ and $c_0+c_1=1$.

\section{ Some properties of the new class of distributions}\label{sec_main}

In this section, we study several probabilistic properties of the density $g(y; \bm{\theta})$ given in \eqref{eq_pdf}.

\subsection{Characteristic function, Mellin transform and moments}
Let $Y$ be a random variable with PDF $g(\cdot; \boldsymbol{\theta})$.
The characteristic function of $Y$, denoted by $\varphi(t)$, is given by
\begin{equation*}
    \varphi(t) = \mathbb{E}_{\boldsymbol{\theta}}\left[e^{itY}\right] = \frac{\mathbb{H}(\theta_1-it, \theta_2, \theta_3, \theta_4, \theta_5, \theta_6)}{\mathbb{H}(\theta_1, \theta_2, \theta_3, \theta_4, \theta_5, \theta_6)}, ~~t\in\mathbb{R}.
\end{equation*}
where $\mathbb{E}_{\boldsymbol{\theta}}(\cdot)$ denotes the expected value under the density $g(\cdot; \boldsymbol{\theta})$.

The Mellin transform\footnote{We may refer to \cite{Springer1979} for further details on the properties of Mellin transforms in the context of the algebra of random variables.} is important in determining the distributions of the products and quotients of independent random variables. Furthermore, in the case of positive random variables, the Mellin transform of the PDF gives us the moments of the random variable.
The Mellin transform of $g(\cdot; \boldsymbol{\theta})$ is given by
$$\{\mathcal{M}g \}(s) = \int_0^\infty y^{s-1} g(y; \boldsymbol{\theta}){\rm d}y = \frac{\mathbb{H}(\theta_1, \theta_2, \theta_3, \theta_4, \theta_5, \theta_6+s-1)}{\mathbb{H}(\theta_1, \theta_2, \theta_3, \theta_4, \theta_5, \theta_6)}.  $$
Then, real moments of $Y$ of order $r$ are given by
\begin{align}\label{moments}
\mathbb{E}_{\boldsymbol{\theta}}[Y^r] 
= 
\{\mathcal{M}g \}(r+1) 
= 
\frac{\mathbb{H}(\theta_1, \theta_2, \theta_3, \theta_4, \theta_5, \theta_6+r)}{\mathbb{H}(\theta_1, \theta_2, \theta_3, \theta_4, \theta_5, \theta_6)}.  
\end{align}

\subsection{Entropy}\label{Entropy}
Let $Y$ be a random variable with PDF $g(\cdot; \boldsymbol{\theta})$. We define the differential entropy of $Y$ by
\begin{equation}\label{eq_entropy1}
    h(Y) :=  -\int_0^\infty g(y; \boldsymbol{\theta}) \log g(y; \boldsymbol{\theta}) {\rm d}y = - \mathbb{E}_{\boldsymbol{\theta}}[\log g(Y; \boldsymbol{\theta})],
\end{equation}
where $\mathbb{E}_{\boldsymbol{\theta}}(\cdot)$ denotes the expected value under the density $g(\cdot; \boldsymbol{\theta})$. Substituting \eqref{eq_pdf} in \eqref{eq_entropy1}, we obtain
\begin{equation*}
    h(Y) = -\int_0^\infty  \left[ -\log c(\boldsymbol{\theta}) + \theta_6 \log y - \theta_1 y - (\theta_2 y^{\theta_3}+\theta_4)^{\theta_5} \right] g(y; \boldsymbol{\theta}) {\rm d}y,
\end{equation*}
which implies
\begin{equation*}
    h(Y) = I_1 + I_2 + I_3 +I_4,
\end{equation*}
where
$$I_1 = \log c(\boldsymbol{\theta}),$$
$$I_2 = -\frac{\theta_6}{c(\boldsymbol{\theta})}\int_0^\infty (\log y)   y^{\theta_6}\exp\left( -\theta_1 y -(\theta_2 y^\theta_3 + \theta_4)^{\theta_5} \right) {\rm d}y,  $$

$$I_3 = \frac{\theta_1}{c(\boldsymbol{\theta})} \int_0^\infty  y^{\theta_6 + 1}\exp\left( -\theta_1 y -(\theta_2 y^\theta_3 + \theta_4)^{\theta_5} \right) {\rm d}y = \theta_1 \frac{\mathbb{H}(\theta_1,\theta_2,\theta_3,\theta_4,\theta_5,\theta_6+1)}{\mathbb{H}(\theta_1,\theta_2,\theta_3,\theta_4,\theta_5,\theta_6)}$$
and 
$$I_4= \frac{1}{c(\boldsymbol{\theta})} \int_0^\infty (\theta_2 y^{\theta_3}+\theta_4)^{\theta_5}  y^{\theta_6}\exp\left( -\theta_1 y -(\theta_2 y^\theta_3 + \theta_4)^{\theta_5} \right) {\rm d}y.$$

It remains to find closed expressions for $I_2$ and $I_4$. For this, consider the case $\theta_5=m\in\mathbb{N}$. Firstly, observe that
$I_2 = -[{\theta_6}/{c(\boldsymbol{\theta})}]\partial_{\theta_6} \mathbb{H}(\boldsymbol{\theta}),$ where $\partial_{\theta_6}$ denotes the partial derivative with respect to $\theta_6$.
It follows from Eq.(16) in \cite{Rathieetal2023} that
\begin{eqnarray}
    \nonumber \partial_{\theta_6} \mathbb{H}(\boldsymbol{\theta}) &=& \frac{1}{\theta_1} \sum_{n=0}^\infty \frac{(-\theta_4^m)^n}{n!} \sum_{k=0}^{mn} \frac{(mn)!}{(mn-k)!k!} \left(\frac{\theta_2}{\theta_4\theta_1^{\theta_3}}\right)^k \partial_{\theta_6}\left\{ \theta_1^{-\theta_6} \Gamma(\theta_6+1+\theta_3 k) \right\}.
\end{eqnarray}
Observe that, for $\theta_1>0$,
\begin{eqnarray}
    \nonumber \partial_{\theta_6}\left\{ \theta_1^{-\theta_6} \Gamma(\theta_6+1+\theta_3 k) \right\} &=& \theta_1^{-\theta_6}\Gamma(\theta_6+1+\theta_3 k) \left[\Psi(\theta_6+1+\theta_3 k) - \log \theta_1  \right], 
\end{eqnarray}
where $\partial_z \Gamma(z) = \Gamma(z)\Psi(z)$. Then,
$$I_2=-\frac{\theta_6\theta_1^{-\theta_6-1}}{c(\boldsymbol{\theta})}\sum_{n=0}^\infty \frac{(-\theta_4^m)^n}{n!} \sum_{k=0}^{mn} \frac{(mn)!}{(mn-k)!k!} \left(\frac{\theta_2}{\theta_4\theta_1^{\theta_3}}\right)^k \Gamma(\theta_6+1+\theta_3 k) \left[\Psi(\theta_6+1+\theta_3 k) - \log \theta_1  \right].$$

Now, for $I_4$, taking $\theta_5=m\in\mathbb{N}$, it follows from Binomial expansion that
\begin{equation*}
I_4 = \frac{1}{c(\boldsymbol{\theta})}\theta_4^m \sum_{k=0}^m \left(\begin{array}{c}
     m  \\
     k 
\end{array}\right) \left(\frac{\theta_2}{\theta_4}\right)^k  \mathbb{H}(\theta_1, \theta_2, \theta_3, \theta_4, m, \theta_3 k +\theta_6).  
\end{equation*}

\subsection{Kullback-Leibler Divergence}\label{KLD-1}

The Kullback-Leibler divergence {(see \cite{Csiszar})} is  useful to measure the difference between two probability distributions.
If $Y_1$ and $Y_2$ are two random variables with PDFs $g_{Y_1}(y;\boldsymbol{\theta})$ and $g_{Y_2}(y;\boldsymbol{\theta}')$, respectively, where $\boldsymbol{\theta}=(\theta_1,\theta_2,\theta_3,\theta_4,\theta_5,\theta_6)$, $\boldsymbol{\theta}'=(\theta_1',\theta_2,\theta_3,\theta_4,\theta_5,\theta_6')$,  $\theta_1\neq\theta_1'$ and $\theta_6\neq\theta_6'$, then their Kullback-Leibler  divergence has the form 
\begin{align*}
D_{\rm KL}(g_{Y_1}\Vert g_{Y_2})
=
\int_0^\infty
g_{Y_1}(y;\boldsymbol{\theta})
\log\left({g_{Y_1}(y;\boldsymbol{\theta})\over g_{Y_2}(y;\boldsymbol{\theta}')}\right)
{\rm d}y.
\end{align*}
From \eqref{eq_pdf} the above integral is
\begin{align*}
    &=\log(c(\boldsymbol{\theta}')-\log(c(\boldsymbol{\theta})
    +
    (\theta_6-\theta_6')\int_0^\infty\log(y)
g_{Y_1}(y;\boldsymbol{\theta}){\rm d}y
+
(\theta_1'-\theta_1)\int_0^\infty y
g_{Y_1}(y;\boldsymbol{\theta}){\rm d}y
\\[0,2cm]
&=\log(c(\boldsymbol{\theta}')-\log(c(\boldsymbol{\theta})
    +
    (\theta_6-\theta_6')
    E_{\boldsymbol{\theta}}[\log(Y)]
+
(\theta_1'-\theta_1)E_{\boldsymbol{\theta}}[Y]
\\[0,2cm]
&=\log(c(\boldsymbol{\theta}')-\log(c(\boldsymbol{\theta})
    +
    (\theta_6-\theta_6')
    \xi_Y
+
(\theta_1'-\theta_1)\, 
\frac{\mathbb{H}(\theta_1, \theta_2, \theta_3, \theta_4, \theta_5, \theta_6+1)}{\mathbb{H}(\theta_1, \theta_2, \theta_3, \theta_4, \theta_5, \theta_6)},
\end{align*}
where $\xi_Y:=E_{\boldsymbol{\theta}}[\log(Y)]$ and in the last equality the formula \eqref{moments} of $E_{\boldsymbol{\theta}}[Y]$ was used. 
That is, 
\begin{align}\label{KLD}
D_{\rm KL}(g_{Y_1}\Vert g_{Y_2})
  =
  \log(c(\boldsymbol{\theta}')-\log(c(\boldsymbol{\theta})
    +
    (\theta_6-\theta_6')
    \xi_Y
+
(\theta_1'-\theta_1)\, 
\frac{\mathbb{H}(\theta_1, \theta_2, \theta_3, \theta_4, \theta_5, \theta_6+1)}{\mathbb{H}(\theta_1, \theta_2, \theta_3, \theta_4, \theta_5, \theta_6)}. 
\end{align}
By using the well-known inequalities $1-x^{-1}<\log(x)\leq x-1$, for $x>0$, from \eqref{moments} we get the finiteness of $\xi_Y$, more precisely,
\begin{align*}
   1-\frac{\mathbb{H}(\theta_1, \theta_2, \theta_3, \theta_4, \theta_5, \theta_6-1)}{\mathbb{H}(\theta_1, \theta_2, \theta_3, \theta_4, \theta_5, \theta_6)} <\xi_Y\leq \frac{\mathbb{H}(\theta_1, \theta_2, \theta_3, \theta_4, \theta_5, \theta_6+1)}{\mathbb{H}(\theta_1, \theta_2, \theta_3, \theta_4, \theta_5, \theta_6)}-1.
\end{align*}

\begin{remark}
   By using the entropy formula (see Subsection \ref{Entropy}) and Kullback-Leibler divergence formula in \eqref{KLD}, we can obtain a closed-form expression for the cross entropy, denoted by $H(g_{Y_1},g_{Y_2})$, since $D_{\rm KL}(g_{Y_1}\Vert g_{Y_2})
=H(g_{Y_1},g_{Y_2})-h(Y_1)$. 
\end{remark}

\section{Estimation}\label{sec_estimation}

In this section, we present two classes of estimators for the model $g(y;\bm\theta)$ \eqref{eq_pdf}.

\subsection{Exponential family}

Consider the case $\theta_5=m\in\mathbb{N}$ being a known parameter and $\theta_3=1$. It follows from Binomial expansion that
\begin{eqnarray}
  \nonumber  g(y; \boldsymbol{\theta}) &=& \exp\left( -\sum_{k=0}^m \left(\begin{array}{c}
         m \\
         k 
    \end{array}\right) \theta_2^k\theta_4^{m-k} y^k - \theta_1 y + \theta_6 \log y - \log c(\boldsymbol{\theta}) \right)\\
    \nonumber &=& \exp\left(\sum_{j=0}^{m+2} a_j(\boldsymbol{\theta}) T_j(y) + b(\boldsymbol{\theta})\right),
\end{eqnarray}
where 
$$a_j(\boldsymbol{\theta}) = \left\{ \begin{array}{cc}
   -\left(\begin{array}{c}
         m \\
         j 
    \end{array}\right) \theta_2^j\theta_4^{m-j},  & j=0,1,\cdots,m,  \\
    -\theta_1, & j=m+1,\\
    \theta_6, &j=m+2,
\end{array} \right.$$
$$T_j(y)=\left\{\begin{array}{cc}
   y^j,  & j=0,1,\cdots,m, \\
    y, & j=m+1,\\
    \log y, & j=m+2, 
\end{array} \right.$$
and $b(\boldsymbol{\theta})=- \log c(\boldsymbol{\theta})$. That means, $g(y; \boldsymbol{\theta})$ belongs to exponential family, provided that $\theta_3=1$ and $\theta_5=m$.

\subsection{Least squares estimation}

{
Let  $\mathbf{Y}=(Y_1, \cdots, Y_n)$ be a random sample of $g(\cdot; \boldsymbol{\theta})$.
Consider the empirical CDF (ECDF) $\hat{G}(y)$ defined as 
$$\hat{G}(y) := \frac{1}{n} \sum_{i=1}^n \mathbbm{1}_{\{Y_i\leq y\}},$$
where $\mathbbm{1}_A$ denotes the indicator of the set $A$.

A widely used estimation method in linear models is the least squares estimator (LSE) (cf. \cite{shao2008mathematical}). In our case, we aim to estimate $\boldsymbol{\theta}$ by comparing the theoretical CDF $G(y; \boldsymbol{\theta})$, given in \eqref{eq_CDF}, with the ECDF $\hat{G}(y)$. To realize this, we minimize the quadratic loss function over a search space $\Theta_0 \subset \Theta$. Thus, the LSE, $\hat{\boldsymbol{\theta}}=\bm\hat{\theta}_{LSE}(\mathbf{Y})$, is given by 
\begin{equation}\label{eq_LSE}
    \hat{\bm\theta} \in \arg\min\left\{ \sum_{j=1}^n \left[G(Y_j; \bm\theta) - \hat{G}(Y_j)\right]^2\right\},
\end{equation}
with respect to $\bm \theta\in\Theta_0$.
}


\subsection{Maximum likelihood estimation}
Let { $\mathbf{Y}=(Y_1, \cdots, Y_n)$}  be a random sample of $g(\cdot; \boldsymbol{\theta})$. The (random) log-likelihood function is given by
\begin{equation}\label{eq_llf}
    \ell(\boldsymbol{\theta}) =\ell(\boldsymbol{\theta}; \mathbf{Y}) = - n\log c(\boldsymbol{\theta}) + \theta_6 \sum_{j=1}^n \log Y_j - \theta_1 \sum_{j=1}^n Y_j  - \sum_{j=1}^n\left(\theta_2 Y_j^{\theta_3} +\theta_4 \right)^{\theta_5}.
\end{equation}

The maximum likelihood estimator (MLE) is given by
\begin{equation}\label{eq_MLE}
    \hat{\boldsymbol{\theta}}_n^{MLE} \in \arg \max_\theta \ell(\boldsymbol{\theta}).
\end{equation}

Observe that \eqref{eq_llf} is a continuous function of $\boldsymbol{\theta}$. Then, on a compact subset $\Theta_0\subset\Theta$, there exists a maximum likelihood estimator.
Computationally, \eqref{eq_llf} is enough to find the MLE.
However, traditional methods of finding MLE via the gradient vector of $\ell(\boldsymbol{\theta})$ can be used. 
{ 
For a more detailed discussion on conditions under which the MLE exists, we refer the reader to \cite{shao2008mathematical}. }

\subsubsection{Partial derivatives of the log-likelihood function}\label{ap_partial derivatives}


Possible candidates for the MLE are the vectors \(\boldsymbol{\theta} \in \Theta\) that satisfy the likelihood equation (cf. \cite{shao2008mathematical})  
\[
\partial_{\boldsymbol{\theta}} \ell(\boldsymbol{\theta}) = \left( \partial_{\theta_1} \ell(\boldsymbol{\theta}),\dots, \partial_{\theta_6} \ell(\boldsymbol{\theta}) \right) = (0,\dots,0).
\]
Since the log-likelihood function \eqref{eq_llf} depends on \( c(\boldsymbol{\theta})=\mathbb{H}(\boldsymbol{\theta}) \), its derivatives, \( \partial_{\boldsymbol{\theta}} c(\boldsymbol{\theta}) \), are needed to explicitly express \( \partial_{\boldsymbol{\theta}} \ell(\boldsymbol{\theta}) \). Moreover, power functions of the form \( x^{\theta_5} \) appear in \eqref{eq_llf}. To simplify the derivatives, we obtain the gradient of the likelihood function when $\theta_5=m\in\mathbb{N}$.
In this case, the MLE can be obtained using numerical procedures for solving the nonlinear system (cf. \cite{burden2010numerical}):
\begin{equation}\label{eq_likelihood_system}
    \partial_{\theta_i} \ell(\boldsymbol{\theta}) = 0, ~~i=1,\cdots,6.
\end{equation}

\begin{remark}
    In the derivatives $\partial_{\boldsymbol{\theta_i}}\ell(\boldsymbol{\theta})$ below, we require $\theta_5 = m \in\mathbb{N}$. In computational terms of applying numerical methods to solve the system \eqref{eq_likelihood_system}, there is no guarantee that $\hat{\theta}_5\in\mathbb{N}$. We could overcome this problem by adding a projection step in the update rule of the numerical algorithm. 
    { That means, we could use the corrected estimator 
    $$\Tilde{\theta}_5 = \arg\inf_{n\in\mathbb{N}} \left| \hat{\theta}_5 - n \right|.$$
    }
\end{remark}


We have the following expressions:
\begin{equation}\label{eq_partial_derivative1_llf}
    \partial_{\theta_1} c(\boldsymbol{\theta}) = - \mathbb{H}(\theta_1, \theta_2, \theta_3, \theta_4, \theta_5, \theta_6-1),
\end{equation}

\begin{equation}\label{eq_partial_derivative2_llf}
    \partial_{\theta_2} c(\boldsymbol{\theta}) = -\theta_5 \sum_{k=1}^{m-1} \left(\begin{array}{c}
         m-1  \\
         k 
    \end{array}\right) \theta_2^k \theta_4^{m-1-k}\mathbb{H}(\theta_1, \theta_2, \theta_3, \theta_4, \theta_5, \theta_6+\theta_3(k+1)), 
\end{equation}

\begin{eqnarray}
    \nonumber     \partial_{\theta_3} c(\boldsymbol{\theta}) &=& -\theta_5\theta_1^{\theta_6+\theta_3-1} \sum_{k=1}^{m-1} \left(\begin{array}{c}
         m-1  \\
         k 
    \end{array}\right) \theta_2^{k+1} \theta_4^{m-1-k}\\
  \label{eq_partial_derivative3_llf}  &\times&\left[ \sum_{l=0}^\infty  \frac{(-\theta_4^m)^l}{l!}\sum_{r=0}^{ml}  \left(\begin{array}{c}
         ml  \\
         r 
    \end{array}\right) \frac{\theta_2^r}{\theta_4^r\theta_1^{3r}} \Gamma(\theta_6+\theta_3+1-\theta_3r)\left( 
\Psi(\theta_6+\theta_3+1-\theta_3r) -\log \theta_1 \right) \right],
\end{eqnarray}

\begin{equation}\label{eq_partial_derivative4_llf}
    \partial_{\theta_4} c(\boldsymbol{\theta}) = - m\sum_{k=0}^{m-1}\left(\begin{array}{c}
         m-1  \\
         k 
    \end{array}\right) \theta_2^{k} \theta_4^{m-1-k}\mathbb{H}(\theta_1, \theta_2, \theta_3, \theta_4, \theta_5, \theta_6+\theta_3k).
\end{equation}

A computable representation of the extreme value $\mathbb{H}$-function, when $\theta_1 > 0$ and $\theta_2 \geq 0$, is 
\begin{equation}\label{ap:eq_approx_H}
    \mathbb{H}(\theta_1, \theta_2, \theta_3, \theta_4, \theta_5,\theta_6) \approx \frac{1}{\theta_1^{1+\theta_6}} \sum_{j=1}^N w_{j,\theta_6} \exp\{-(\theta_2 \theta_1^{-\theta_3} x_j^{\theta_3}+\theta_4)^{\theta_5} \},
\end{equation}
where $x_j$ is the $j$-th root of the generalized Laguerre polynomials $L_N^{(\alpha)} (x)$ and the weight $w_{j,\alpha}$ is given by \cite{abramowitzstegun,Stroud}:
\begin{equation}
    w_{j,\alpha} = \frac{\Gamma(N+\alpha+1) x_j}{N!(N+1)^2 [L_{N+1}^{(\alpha)}(x_j)]^2}.
\end{equation}
This way, the derivative with respect to $\theta_5$ can be given as:
\begin{eqnarray}\label{deriv}
    \partial_{\theta_5} c(\boldsymbol{\theta}) \approx  - \frac{1}{\theta_1^{1+\theta_6}} \sum_{j=1}^N \frac{w_{j,\theta_6} (\theta_2 \theta_1^{-\theta_3} x_j^{\theta_3}+\theta_4)^{\theta_5} \text{log} \left(\theta_2 \theta_1^{-\theta_3} x_j^{\theta_3}+\theta_4 \right)}{\exp\{(\theta_2 \theta_1^{-\theta_3} x_j^{\theta_3}+\theta_4)^{\theta_5} \}}.
\end{eqnarray}
In particular, if $\theta_5=m\in\mathbb{N}$, 
\begin{equation}\label{eq_partial_derivative5_llf}
    \partial_{\theta_5} c(\boldsymbol{\theta}) \approx - \frac{1}{\theta_1^{1+\theta_6}} \sum_{j=1}^N \sum_{k=0}^m \left(\begin{array}{c}
         m  \\
         k 
    \end{array}\right) \theta_2^k \theta_1^{-\theta_3 k} \theta_4^{m-k} x_j^k \frac{w_{j,\theta_6} \text{log} \left(\theta_2 \theta_1^{-\theta_3} x_j^{\theta_3}+\theta_4 \right)}{\exp\{(\theta_2 \theta_1^{-\theta_3} x_j^{\theta_3}+\theta_4)^{m} \}}.
\end{equation}
Finally,
\begin{equation}\label{eq_partial_derivative6_llf}
    \partial_{\theta_6} c(\boldsymbol{\theta}) = \theta_1^{-\theta_6-1}\sum_{n=0}^\infty \frac{(-\theta_4^m)^n}{n!} \sum_{k=0}^{mn} \frac{(mn)!}{(mn-k)!k!} \left(\frac{\theta_2}{\theta_4\theta_1^{\theta_3}}\right)^k \Gamma(\theta_6+1+\theta_3 k) \left[\Psi(\theta_6+1+\theta_3 k) - \log \theta_1  \right].
\end{equation}

Now, we are able to describe the score functions:
\begin{equation*}
    \partial_{\theta_1} \ell(\boldsymbol{\theta}) = n \frac{\mathbb{H}(\theta_1, \theta_2, \theta_3, \theta_4, \theta_5, \theta_6 - 1)}{\mathbb{H}(\theta_1, \theta_2, \theta_3, \theta_4, \theta_5, \theta_6)} - \sum_{j=1}^n Y_j,
\end{equation*}
\begin{eqnarray}
   \nonumber \partial_{\theta_2} \ell(\boldsymbol{\theta}) &=& - \frac{n}{c(\boldsymbol{\theta})}\left[ -\theta_5 \sum_{k=0}^{m-1}\left(\begin{array}{c}
         m-1 \\
         k
    \end{array}\right)\theta_2^k\theta_4^{m-1-k} \mathbb{H}(\theta_1, \theta_2, \theta_3, \theta_4, \theta_5, \theta_6 + \theta_3(k+1))\right] \\
  \nonumber  &-&\theta_5 \sum_{j=1}^n \left(\theta_2 Y_j^{\theta_3} + \theta_4\right)^{\theta_5 -1} Y_j^{\theta_3},
\end{eqnarray}
\begin{equation*}
     \partial_{\theta_3} \ell(\boldsymbol{\theta}) = - \frac{n}{c(\boldsymbol{\theta})} \partial_{\theta_3} c(\boldsymbol{\theta}) - m\theta_2 \sum_{j=1}^n\sum_{k=0}^{m-1}\left(\begin{array}{c}
          m-1 \\
          k 
     \end{array}\right) \theta_2^k \theta_4^{m-1-k} Y_j^{\theta_3(k+1)} \log Y_j,
\end{equation*}
where $\partial_{\theta_3} c(\boldsymbol{\theta})$ is given in \eqref{eq_partial_derivative3_llf}.
\begin{equation*}
     \partial_{\theta_4} \ell(\boldsymbol{\theta})= - m\sum_{k=0}^{m-1}\left(\begin{array}{c}
          m-1  \\
          k
     \end{array}\right) \theta_2^k\theta_4^{m-1-k} \mathbb{H}(\theta_1, \theta_2, \theta_3, \theta_4, \theta_5, \theta_6 + \theta_3k),
\end{equation*}
\begin{eqnarray}
     \nonumber \partial_{\theta_5} \ell(\boldsymbol{\theta}) &\approx& - \frac{n}{c(\boldsymbol{\theta})}\left[\frac{1}{\theta_1^{1+\theta_6}} \sum_{j=1}^N \sum_{k=0}^m \left(\begin{array}{c}
         m  \\
         k 
    \end{array}\right) \theta_2^k \theta_1^{-\theta_3 k} \theta_4^{m-k} x_j^k \frac{w_{j,\theta_6} \text{log} \left(\theta_2 \theta_1^{-\theta_3} x_j^{\theta_3}+\theta_4 \right)}{\exp\{(\theta_2 \theta_1^{-\theta_3} x_j^{\theta_3}+\theta_4)^{m} \}}\right]\\
    \nonumber &-& \sum_{j=1}^n \sum_{k=0}^m \left(\begin{array}{c}
         m  \\
         k 
    \end{array}\right) \theta_2^k\theta_4^{m-k} Y_j^{\theta_3k}
    \log(\theta_2 Y_j^{\theta_3}+\theta_4),
\end{eqnarray}
and 
\begin{equation*}
     \partial_{\theta_6} \ell(\boldsymbol{\theta}) = -\frac{n}{c(\boldsymbol{\theta})} \partial_{\theta_6} c(\boldsymbol{\theta}) + \sum_{j=1}^n\log Y_j,
\end{equation*}
where $\partial_{\theta_6} c(\boldsymbol{\theta})$ is given in \eqref{eq_partial_derivative6_llf}.

\begin{remark}
    The gradient $\nabla \ell(\boldsymbol{\theta})$ can be {computed} for the general case $\theta_5\in\mathbb{R}$ using approximation \eqref{ap:eq_approx_H}. 
\end{remark}

\subsection{Computational implementation}

\subsubsection{{Choosing a good initial guess}}
{ While the methods used for estimation may be straightforward in the sense that they rely on standard estimation techniques, the presence of integral expressions in the PDF means that the computational implementation/optimization of models involving special functions is not always simple. Optimizers in languages such as Python~\cite{python32009} and R~\cite{Rsoftware} can easily become inefficient if a well-planned strategy for implementing objective functions is not adopted.
A good initial guess will be important for efficient convergence of the implemented methods.

\begin{ex}  \label{ex:initial_guess}
    Let $\mathbf{Y} = (Y_1, \dots, Y_n)$ be a sample from the distribution \( g(\cdot; \bm\theta) \). Suppose we estimate \( (\hat{\alpha}, \hat{\beta}) \) by modeling \( \mathbf{Y} \) as a Gamma\((\alpha, \beta)\) distribution and find that this model fits the data well. Then, from Table \ref{tab:g_particular}, we obtain the initial estimate \( \hat{\boldsymbol{\theta}}_0 = (\hat{\beta}, 0, 1, 1, \hat{\alpha} - 1) \) for the parameter vector \( \boldsymbol{\theta} \) in the general model \( g(\cdot; \bm\theta) \).  

    In the next step, we can use \( \hat{\boldsymbol{\theta}}_0 \) as an initial guess to refine our estimation of \( \boldsymbol{\theta} \) without being constrained to the specific Gamma model. This procedure can also be applied using other distributions listed in Table \ref{tab:g_particular}.  
\end{ex}  

Although simple, Example \ref{ex:initial_guess} provides a powerful tool for obtaining initial guesses, as the distributions in Table \ref{tab:g_particular} are mostly available in standard Python~\cite{python32009} and R~\cite{Rsoftware} libraries.
Another effective approach for obtaining an initial guess is to use the result of one estimation method (such as LSE, method of moments, methods based on characteristic function or Mellin transform) as the starting point for another method (such as MLE). This will be discussed in the next subsection.

Finally, a different option for choosing the initial guess involves defining a grid of values and testing various initial values within a reasonable range. The estimation method's performance is then evaluated for each of these values. The aim is to identify the value from the grid that leads to the most accurate final estimate, such as the one with the best convergence or smallest error. This process can be carried out empirically by observing which initial guess produces the most favorable result.
}

\subsubsection{ Algorithms for estimations}

The parameter estimation was done by optimization procedures according to Equation~\eqref{eq_LSE} for the LSE and to Equation~\eqref{eq_MLE} for MLE.

The large number of parameters demanded by extreme value { $\mathbb{H}$-function} class of distributions lead to frequent convergence problems when the MLE approach was tried even using particular cases of Table~\ref{tab:g_particular} as initial guesses.

Thus, we decided to use LSE in Algorithm~\ref{algo_lse} as a first step in the parameter estimation due to its better convergence and try to refine the estimation by using a LSE-tuple as an initial guess to the MLE estimator of Algorithm~\ref{algo_mle}. 

\begin{algorithm}
\caption{Estimation of parameters $\bm\theta$ using the LSE}
\label{algo_lse}
\begin{algorithmic}[1]
    \Statex \textbf{Input:} Data $(Y_1, \cdots, Y_n)$ sampled from the $g(\cdot; \bm\theta)$.
    \Statex \textbf{Output:} Estimates $\widehat{\theta_{Y}}$.
    
    \State Compute estimates $\hat{\alpha}$, $\hat{\beta}$, $\hat{\gamma}$, and $\hat{\sigma}$ to fit the particular case chosen from Table~\ref{tab:g_particular}.
    
    \State Set the initial guess for the numerical optimization to be $\bm\theta_0 = (\theta_1,\theta_2,\theta_3,\theta_4,\theta_5,\theta_6)$ also according to Table~\ref{tab:g_particular}.
    
    \State Calculate $\widehat{\theta_{Y}}$ using \eqref{eq_LSE}.
    
    \State \textbf{Return} $\widehat{\theta_{Y}}$.
\end{algorithmic}
\end{algorithm}




\begin{algorithm}
\caption{Estimation of parameters $\bm\theta$ using the MLE}
\label{algo_mle}
\begin{algorithmic}[1]
    \Statex \textbf{Input:} Data $(Y_1, \cdots, Y_n)$ sampled from the $g(\cdot; \bm\theta)$.
    \Statex \textbf{Output:} Estimates $\hat{\theta}^{MLE}$.
    
    \State Set the initial guess for the numerical optimization to be $\bm\theta_0 = \widehat{\theta_{Y}}$ from Algorithm \ref{algo_lse}.
    
    \State Calculate $\hat{\theta}^{MLE}$ using \eqref{eq_MLE}.
    
    \State \textbf{Return} $\hat{\theta}^{MLE}$.
\end{algorithmic}
\end{algorithm}



  In the next section, we apply Algorithms~\ref{algo_lse} and~\ref{algo_mle} to model real data.
  The algorithms were implemented using the Python language~\cite{python32009} and the optimizers \verb|optimize.scipy|~\cite{2020SciPy-NMeth} and \verb|hyperopt|~\cite{bergstra2013making}.
 Furthermore, the R software \cite{Rsoftware} was used to obtain descriptive statistics and the figures presented throughout this article (with the exception of Figure~\ref{fig:fits_extreme_h}, which was generated in Python). All the code and a tutorial on how to use it are available at \url{https://github.com/eip-unb/positive_support_asymmetric}
.

\section{Applications}\label{sec_applications}

In order to show the versatility of extreme value $\mathbb{H}$-function class of distributions, 
{ we show that the model was able to fit three datasets with distinct characteristics. 
Until now, these datasets had been effectively modeled by three different and competing models \cite{lawless2011statistical, Ramos2020Frechet, valiollahi2013estimation}. Here, we unify them into the model \eqref{eq_CDF}. Even when accounting for penalties due to the increased number of parameters, the new model is at least as effective as its particular cases, demonstrating the versatility of the model.}

The study by \cite{Ramos2020Frechet} shown that Frechét distributions fitted better the minimum monthly flows of water ($m^3/s$) on the Piracicaba Riber, located in São Paulo state, Brazil, from 1960 to 2014. We took the minimum monthly flow for September as presented below:

\begin{eqnarray}
  \nonumber  X&=&( 29.19, 8.49, 7.37, 82.93, 44.18, 13.82, 22.28, 28.06, 6.84, 12.14, \\
    \nonumber && 153.78, 17.04, 13.47, 15.43, 30.36, 6.91, 22.12, 35.45, 44.66, 95.81, \\
    \nonumber && 6.18, 10.00, 58.39, 24.05, 17.03, 38.65, 47.17, 27.99, 11.84, 9.60, \\
    \nonumber && 6.72, 13.74, 14.60, 9.65, 10.39, 60.14, 15.51, 14.69, 16.44).
\end{eqnarray}

The strength of carbon fibers tested under stressing tensions is frequently used in the literature \cite{valiollahi2013estimation} and we take the strength data measured in GPa (Gigapascal) for single carbon fibers of 20 mm length as presented below
\begin{eqnarray}
  \nonumber  Y&=& (1.312, 1.314, 1.479, 1.552, 1.700, 1.803, 1.861, 1.865, 1.944, 1.958, 1.966, \\
  \nonumber &&1.977, 2.006, 2.021, 2.027, 2.055, 2.063, 2.098, 2.140, 2.179, 2.224, 2.240, \\
  \nonumber && 2.253, 2.270, 2.272, 2.274, 2.301, 2.301, 2.359, 2.382, 2.382, 2.426, 2.434, \\
  \nonumber && 2.435, 2.478, 2.490, 2.511, 2.514, 2.535, 2.554, 2.566, 2.570, 2.586, 2.629, \\
  \nonumber && 2.633, 2.642, 2.648, 2.684, 2.697, 2.726, 2.770, 2.773, 2.800, 2.809, 2.818, \\
  \nonumber && 2.821, 2.848, 2.880, 2.954, 3.012, 3.067, 3.084, 3.090, 3.096, 3.128, 3.233,\\
  \nonumber && 3.433, 3.585, 3.585)
\end{eqnarray}
to be modelled as a Weibull RV.

Life data are sometimes modeled with the Gamma distribution and \cite{lawless2011statistical} or \cite{nelson2005applied} discuss applications of the gamma distribution to life data.
The following dataset represents failure times of machine parts, some of which are manufactured by two manufacturers and were merged after a previous test which rejected the hypothesis of significant differences in the part life for the two manufacturers.

\begin{eqnarray}
  \nonumber  Z&=& ( 620, 470, 260, 89, 388, 242, 103, 100, 39, 460, 284, 1285, 218, 393, 106, \\
  \nonumber && 158, 152, 477, 403, 103, 69, 158, 818, 947, 399, 1274, 32, 12, 134, 660, 548, \\
  \nonumber && 381, 203, 871, 193, 531, 317, 85, 1410, 250, 41, 1101, 32, 421, 32, 343, 376, \\
  \nonumber && 1512, 1792, 47, 95, 76, 515, 72, 1585, 253, 6, 860, 89, 1055, 537, 101, 385,  \\
  \nonumber && 176, 11, 565, 164, 16, 1267, 352, 160, 195, 1279, 356, 751, 500, 803, 560, 151, \\
  \nonumber && 24, 689, 1119, 1733, 2194, 763, 555, 14, 45, 776, 1, 1747, 945, 12, 1453, 14, \\
  \nonumber && 150, 20, 41, 35, 69, 195, 89, 1090, 1868, 294, 96, 618, 44, 142, 892, 1307, 310, \\
  \nonumber && 230, 30, 403, 860, 23, 406, 1054, 1935, 561, 348, 130, 13, 230, 250, 317, 304, \\
  \nonumber && 79, 1793, 536, 12, 9, 256, 201, 733, 510, 660, 122, 27, 273, 1231, 182, 289, \\
  \nonumber && 667, 761, 1096, 43, 44, 87, 405, 998, 1409, 61, 278, 407, 113, 25, 940, 28, \\
  \nonumber && 848, 41, 646, 575, 219, 303, 304, 38, 195, 1061, 174, 377, 388, 10, 246, 323, \\
  \nonumber && 198, 234, 39, 308, 55, 729, 813, 1216, 1618, 539, 6, 1566, 459, 946, 764, \\
  \nonumber && 794, 35, 181, 147, 116, 141, 19, 380, 609, 546).
\end{eqnarray}

Descriptive statistics for the $X, Y$ and $Z$ RVs are presented in Table~\ref{tab:datasets_descritive}. For each dataset, the following summary statistics are given: the minimum value observed in the dataset (Min.);
the first quartile value (1st Qu.); the median value (Median); the average value (Mean); the third quartile value (3rd Qu.); the maximum value observed (Max.); the standard deviation (Sd.); the coefficient of Symmetry (CS); the coefficient of Kurtosis (CK);  the sample size (n).


\begin{table}[H]
\caption{{ Descriptive statistics for the RVs $X, Y$ and $Z$.}}
\label{tab:datasets_descritive}
\begin{tabular}{rrrrrrrrrrr}
\hline
Dataset & Min.  & 1st Qu. & Median & Mean       & 3rd Qu. & Max.   & Sd         & CS          & CK          & n   \\
\hline
$X$       & 6.180  & 11.115  & 16.440  & 28.285 & 32.905  & 153.78 & 29.319 & 2.451  & 6.785  & 39  \\
$Y$       & 1.312 & 2.098   & 2.478  & 2.451 & 2.773   & 3.585  & 0.495 & -0.027 & -0.148 & 69  \\
$Z$       & 1     & 96      & 304    & 463.647 & 667     & 2194   & 476.956 & 1.358  & 1.250  & 201\\
\hline
\end{tabular}
\end{table}

\begin{figure}[H]
    \centering
    \subfloat[]{\includegraphics[width=.52\linewidth]{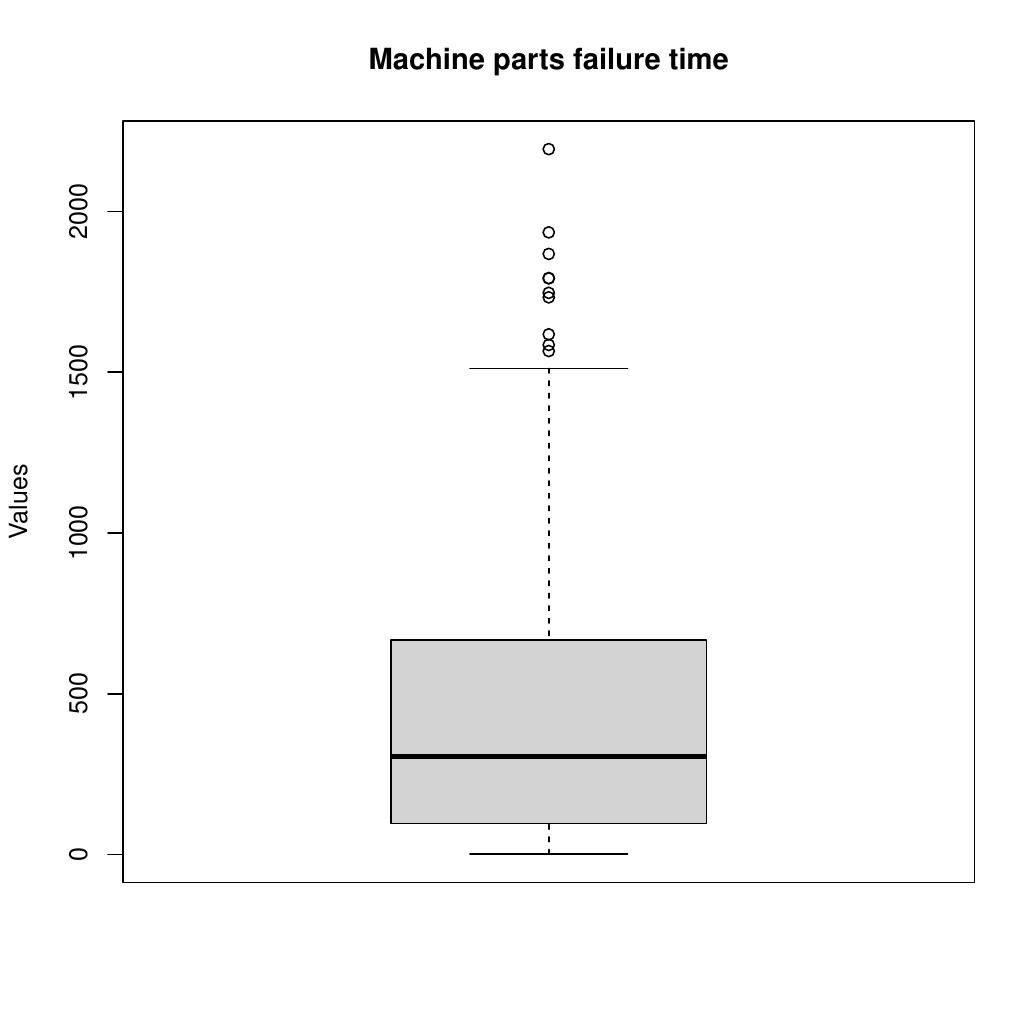}}
    \subfloat[]{\includegraphics[width=.52\linewidth]{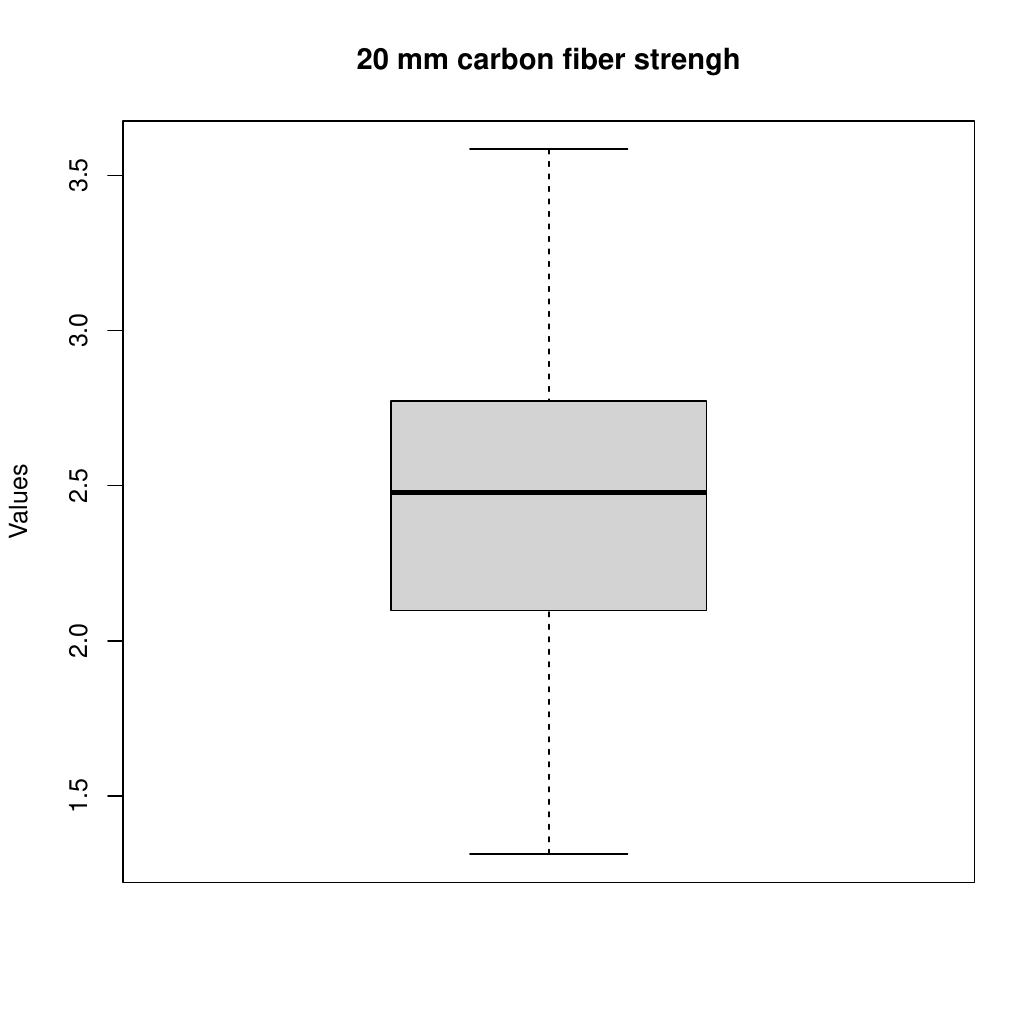}} \\
    \subfloat[]{\includegraphics[width=.52\linewidth]{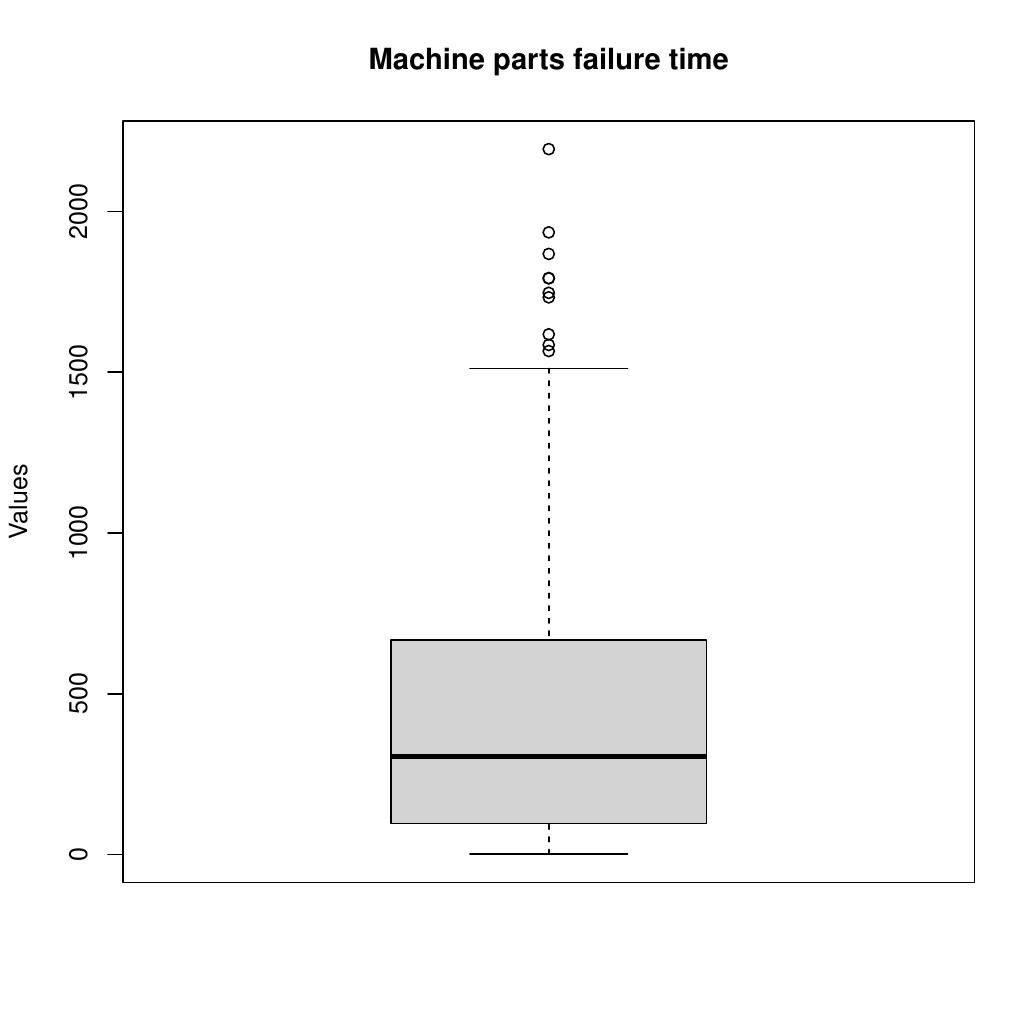}}
    \caption{Boxplot for (a) $X$, (b) $Y$ and (c) $Z$.}
    \label{fig:boxplot}
\end{figure}

    The descriptive statistics (Table \ref{tab:datasets_descritive}) and boxplots (Figure \ref{fig:boxplot}) show that the distributions of data $X$ and $Z$ are right-skewed. This supports the initial choice of using the Fréchet, Weibull, and Gamma models for each case, respectively.

\begin{figure}[H]
    \centering
    \subfloat[]{\includegraphics[width=.52\linewidth]{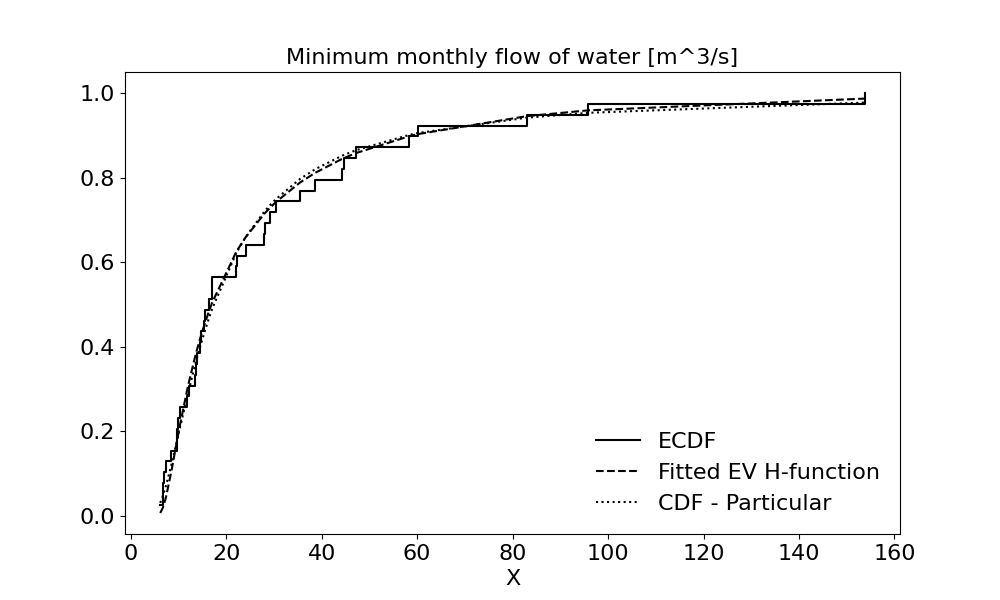}}
    \subfloat[]{\includegraphics[width=.52\linewidth]{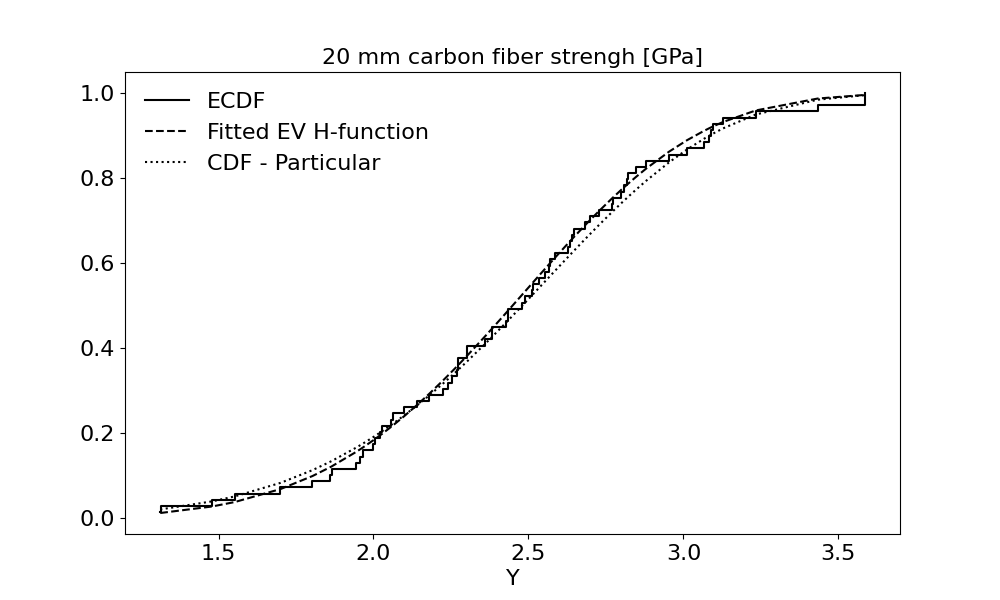}} \\
    \subfloat[]{\includegraphics[width=.52\linewidth]{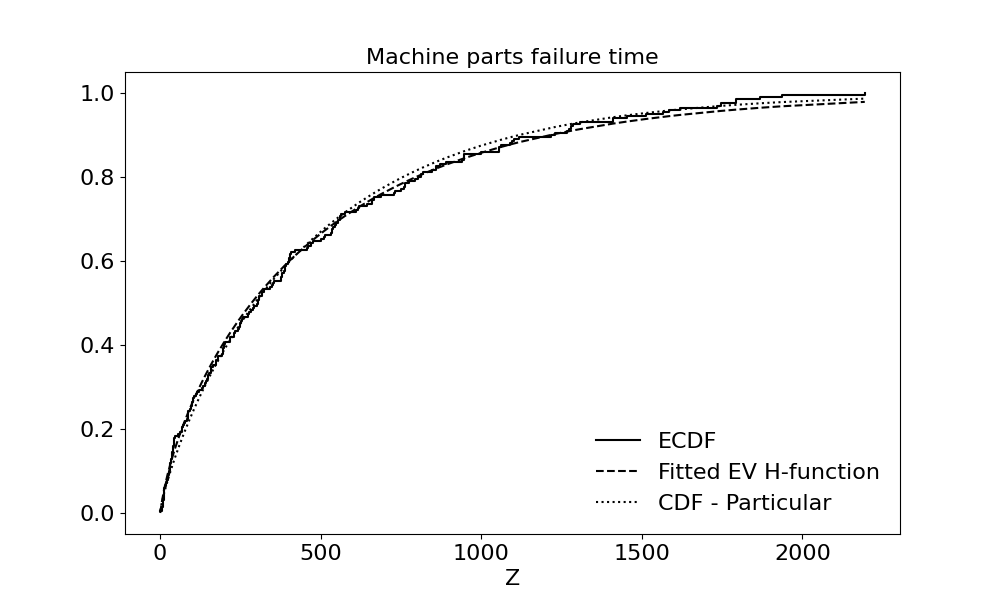}}
    \caption{Empirical and fitted extreme value $\mathbb{H}$-function distribution functions for (a) $X$, (b) $Y$ and (c) $Z$.}
    \label{fig:fits_extreme_h}
\end{figure}

We aim to verify that the newly developed model in this study not only generalizes the existing models mentioned earlier but also provides a better fit to the data. 
Figure~\ref{fig:fits_extreme_h} shows the fit of distributions to dataset, meanwhile the estimated parameters of the previous modeling reported by \cite{Ramos2020Frechet, valiollahi2013estimation} were $(\hat{\alpha}, \hat{\sigma})$ equal to (1.564, 13.76) for X (Fréchet model), (5.505, 2.651) for Y (Weibull model) and (0.8274, 560.35) for Z (Gamma model).

According to Table \ref{tab:g_particular}, the Fréchet RV $X$ can be expressed by a distribution of the extreme value { $\mathbb{H}$-function} class with the parameter set $\bm\theta_{X_0} = (0, 0.0727, 1, 0, -1.564, -2.564)$. 
We take $\bm\theta_{X_0}$ as a initial guess and estimate the parameter set $\widehat{\bm\theta_{X}}$ which better fits to data.
The parameters estimation was done by means of an optimization procedure by minimizing the distance (mean squared error) between the empirical cumulative density function (ECDF) and the CDF (see Algorithm \ref{algo_lse}).


{The parameter set $\hat{\theta}^{MLE}_X = (0.021, 0.0008, 3.973, 0.007, -9.51, -0.899)$ yields a better fit as supported by the information criteria such as the Akaike Information Criterion (AIC), Bayesian Information Criterion (BIC), and Efficient Determination Criterion (EDC). 
These estimates were obtained using Algorithm \ref{algo_mle}.
The same procedure were used to get $\widehat{\theta_{Y}}$ and $\widehat{\theta_{Z}}$ (respectively, $\hat{\theta}^{MLE}_Y$ and $\hat{\theta}^{MLE}_Z$),  which are presented on Table~\ref{tab:tab_estimation_results}.}

   

\begin{table}[H]
\caption{{Estimated parameters and information criteria for model selection. $\mathbb{H}$-extreme distribution outperforms particular cases. } }
\label{tab:tab_estimation_results}
\centering
\begin{tabular}{cc|c|ccc}
  \hline
Dataset & PDF & $\widehat{\theta}$ &  AIC & BIC & EDC \\  
  \hline
  $X$  & Fréchet  & $\bm\theta_{X_0} = (0, 0.0727, 1, 0, -1.564, -2.564)$ & 333.70 & 343.68 & 313.9 \\
  & \textbf{$\mathbb{H}$-extreme}  & $\widehat{\bm\theta_{X}} = (0.013, 0.364, 0.517, 0.0006, -5.089, -1.446)$ & 332.34 & 342.32 & 312.55 \\
    &   & $\hat{\theta}^{MLE}_X = (0.021, 0.0008, 3.973, 0.007, -9.51, -0.899)$ & \textbf{328.89} & \textbf{338.87} & \textbf{309.10} \\


    
  \hline
     $Y$ &  Weibull  & $\bm\theta_{Y_0} = (0, 0.377, 1, 0, 5.50, 4.50)$ & 111.19 & 124.60 & 90.53 \\
    & \textbf{$\mathbb{H}$-extreme}  & $\widehat{\bm\theta_{Y}} = (0.002, 0.5, 0.80, 0.03, 5.44, 6.12)$ & 110.35 & 123.75 & 89.69 \\
        &  & $\hat{\bm\theta}^{MLE}_X = (0.0, 0.657, 0.621, 0.0,
 5.394, 7.150)$ & \textbf{109.75} & \textbf{123.15} & \textbf{89.09} \\
   \hline
    $Z$ & Gamma  & $\bm\theta_{Z_0} = (0.0018, 0, 1, 0, 1, -0.1726)$ & 2876.85 & 2896.67 & 2854.84 \\
    & \textbf{$\mathbb{H}$-extreme}  & $\widehat{\bm\theta_{Z}} = (0.0015, 7.17, -6.55, 0.0018, 8.57, -0.29)$ & 2865.10 & 2884.92 & 2843.09 \\
        &   & $\hat{\bm\theta}^{MLE}_Z = (0.0017, 7.17, -6.55, 0.005, 8.57, -0.195)$ & \textbf{2863.10} & \textbf{2882.92} & \textbf{2841.09} \\

        

        
   \hline
   
\end{tabular}
\end{table}
 
\begin{figure}[htb!]
    \centering
    \subfloat[]{\includegraphics[width=.52\linewidth]{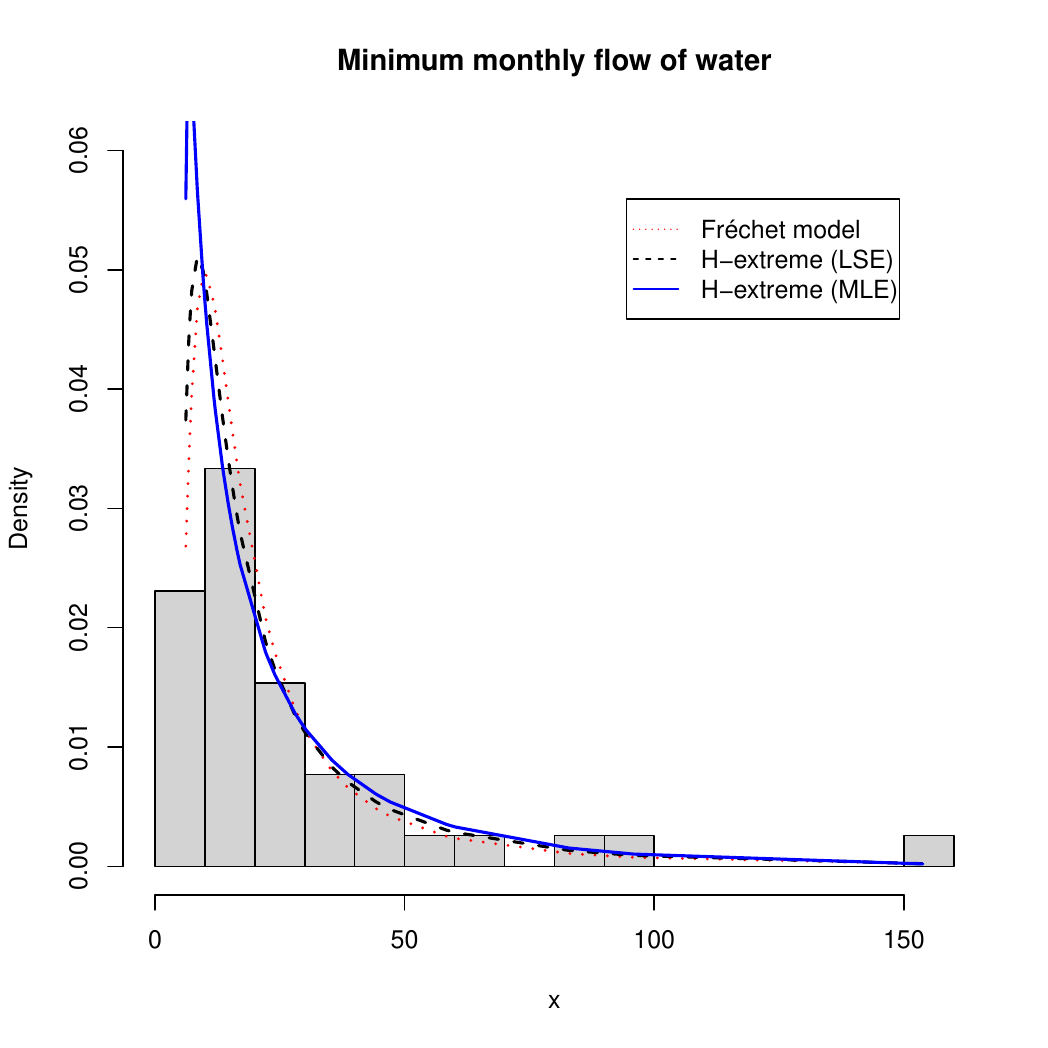}}
    \subfloat[]{\includegraphics[width=.52\linewidth]{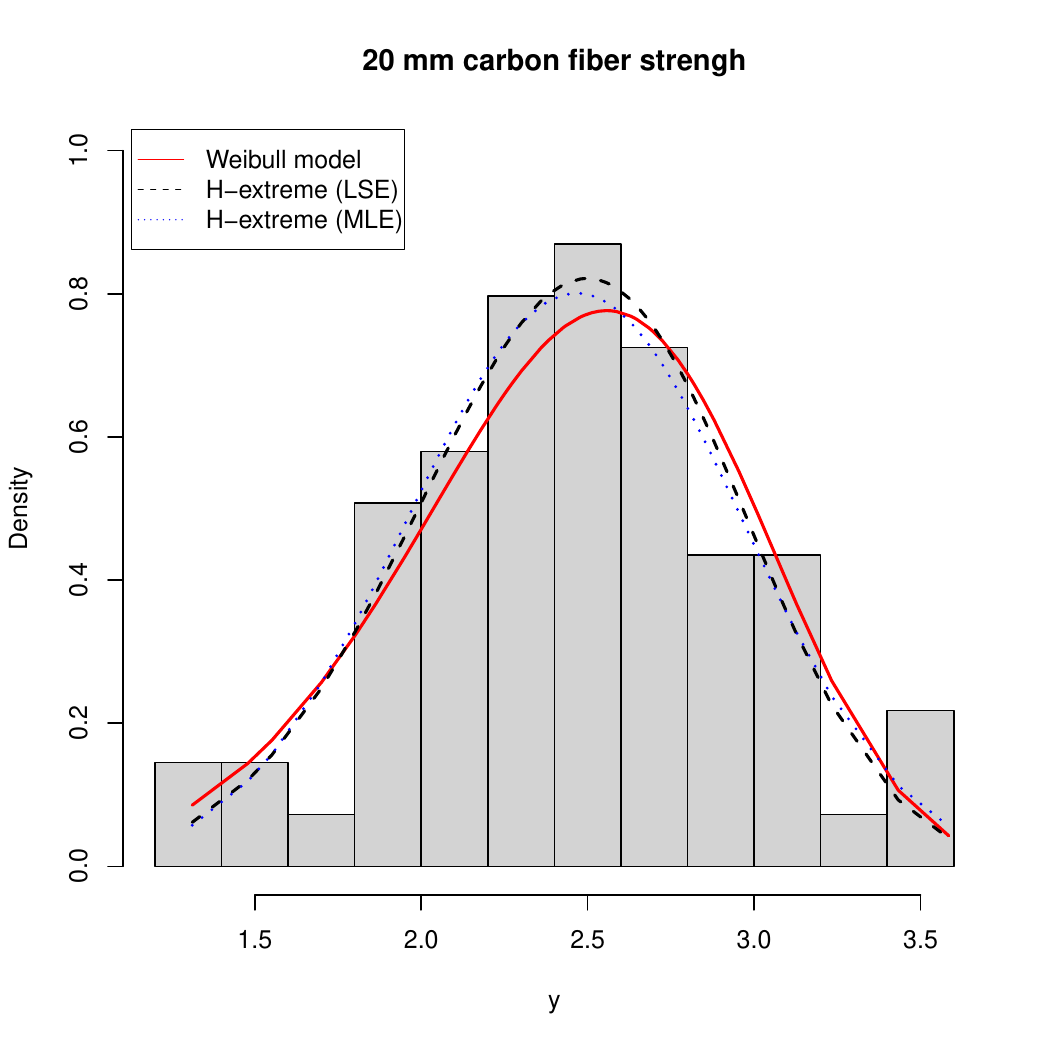}} \\
    \subfloat[]{\includegraphics[width=.52\linewidth]{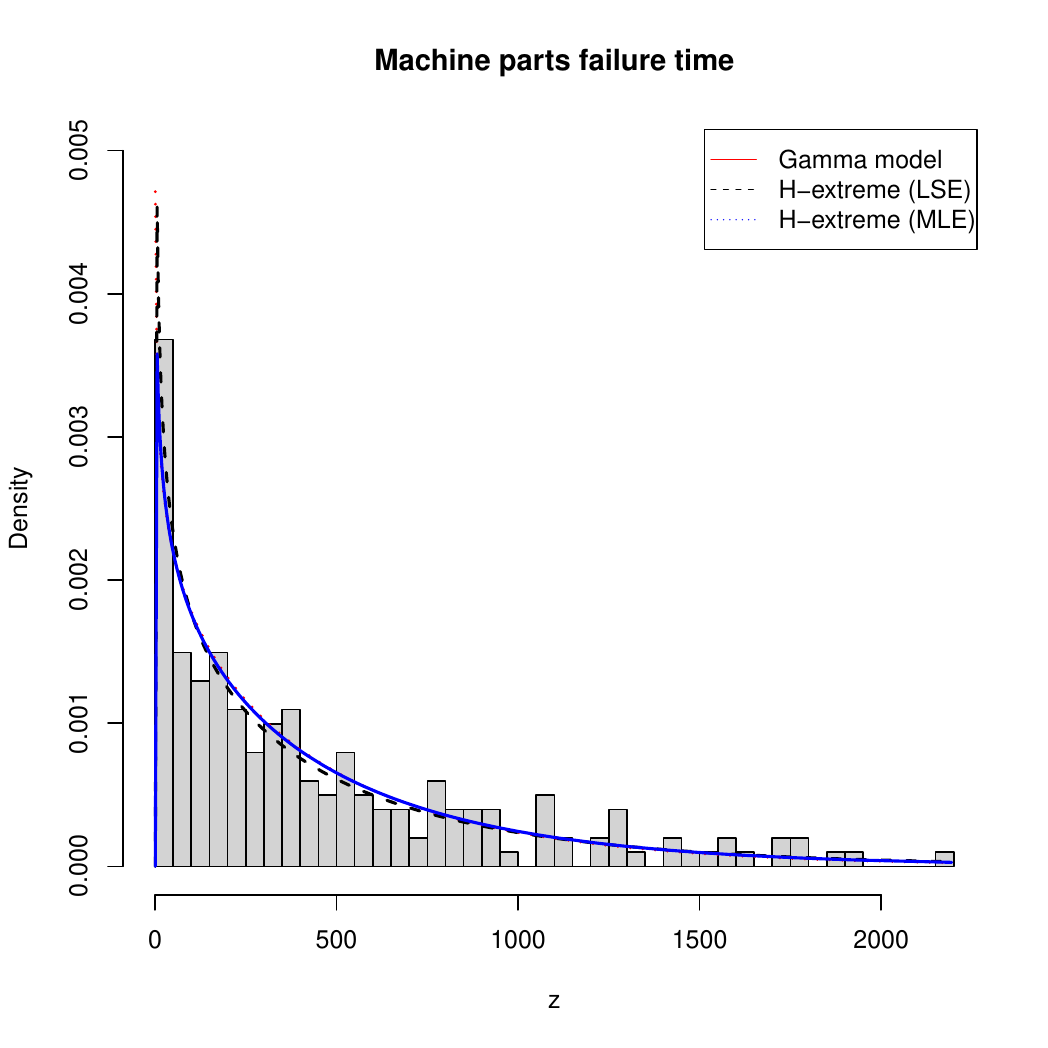}}
    \caption{Histograms and fitted probability density functions for (a) $X$, (b) $Y$ and (c) $Z$.}
    \label{fig:hist}
\end{figure}

{Table \ref{tab:tab_estimation_results} summarizes the improvements achieved through estimations using both LSE (Algorithm \ref{algo_lse}) and MLE (Algorithm \ref{algo_mle}). For all the datasets, the modeling gains of using MLE instead of LSE were validated by the information criteria (AIC, BIC, EDC). Such findings are justified by the fact that MLE explicitly incorporates the probabilistic structure of the data through the likelihood function. By maximizing the likelihood, MLE identifies parameter values that make the observed data most probable under the assumed model, leading to estimators that are typically unbiased, consistent, and asymptotically efficient. In contrast, LSE focuses solely on minimizing the squared deviation between observed and predicted values, without considering the underlying distributional properties. As a result, MLE tends to provide more accurate and statistically meaningful estimates. Figure \ref{fig:hist} shows the fit of the theoretical PDF to the histogram of the data.}

To ensure statistical validity when model parameters are estimated from the data, a parametric bootstrap procedure was implemented to obtain reliable p-values for the fitted CDFs. Specifically, we employed the Kolmogorov–Smirnov (KS) and Cramér–von Mises (CVM) goodness-of-fit tests, combined with bootstrap-based resampling using {
$M=1,000$} (see, e.g., Chapter 4.2.3 in \cite{davison1997bootstrap}). This approach accounts for the composite nature of the null hypothesis and provides corrected significance levels. The resulting p-values indicate that the fitted models adequately describe the data: {for dataset $X$, the KS and CVM p-values are 0.60 and 0.83, respectively; for $Y$, 0.99 and 0.99; and for $Z$, 0.99 and 0.95.} All values are sufficiently high to suggest no significant departure between the empirical and fitted distributions.



\begin{figure}[H]
    \centering
    \subfloat[]{\includegraphics[width=0.9\linewidth, height=0.3\textheight]{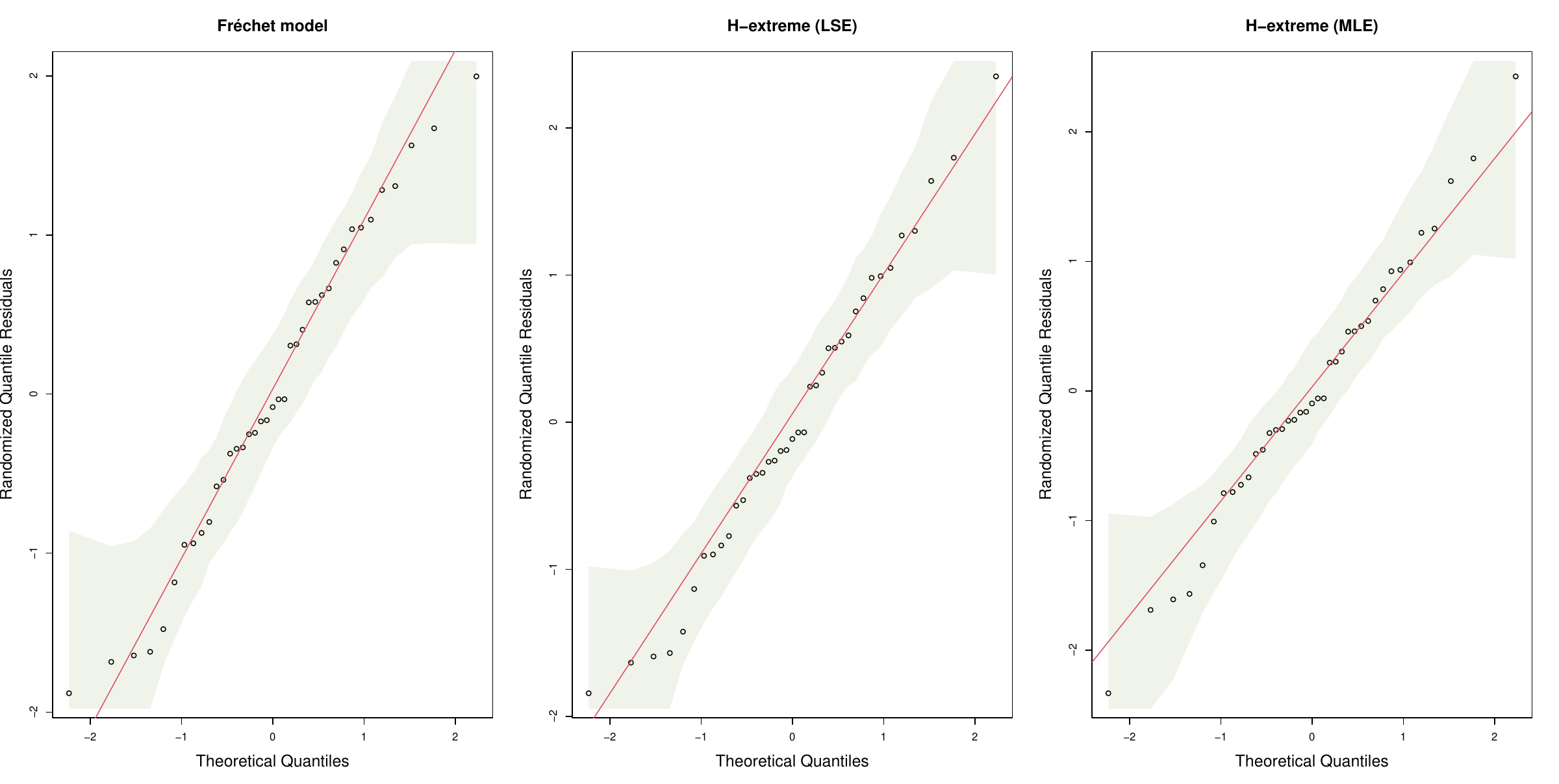}}\\
    \subfloat[]{\includegraphics[width=0.9\linewidth, height=0.3\textheight]{Figs/qq_x.pdf}} \\
    \subfloat[]
    {\includegraphics[width=0.9\linewidth, height=0.3\textheight]{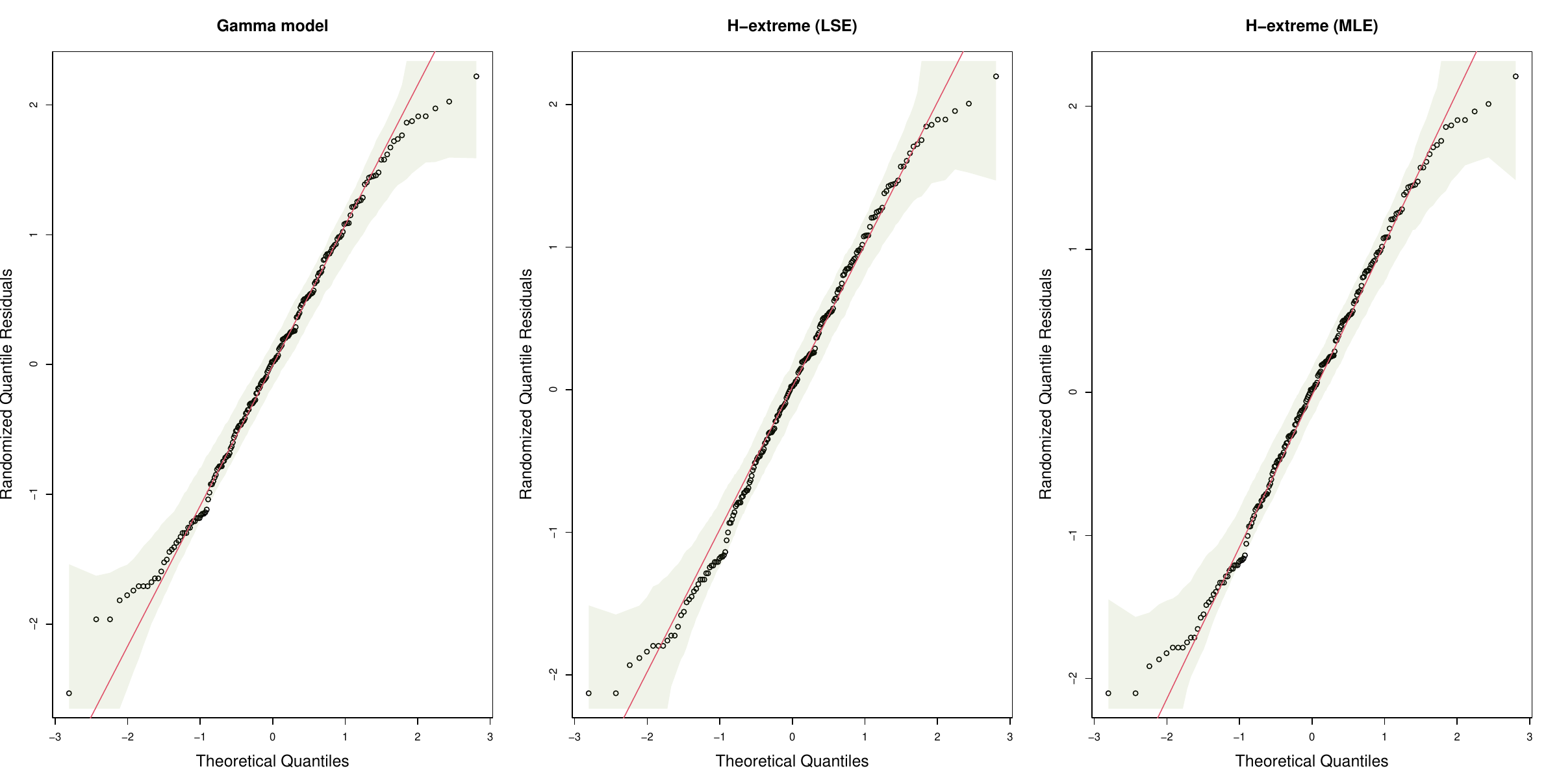}}
    \caption{Normal Quantile-Quantile plot displaying residuals from fitted models for (a) $X$, (b) $Y$ and (c) $Z$.}
    \label{fig:qq}
\end{figure}

To ensure a good model fit, we use Randomized Quantile (RQ) residuals, as defined by \cite{dunn1996randomized}. The RQ residuals are calculated using the formula 
$R_i=\Phi^{-1}(G(y_i; \hat{\bm \theta}))$, where $\hat{\bm \theta}$ is the estimated parameter, $G(y_i; \hat{\bm \theta})$ is the CDF \eqref{eq_CDF} of the model fitted to each observation \(y_i\), and $\Phi^{-1}$ represents the quantile of the standard normal distribution $N(0, 1)$. When $F$ is continuous, the RQ residuals follow a standard normal distribution, excluding the impact of sampling variability in the estimated parameters.





{For all datasets, the Quantile–Quantile (QQ) plots indicate that the $\mathbb{H}$-extreme distribution provides a superior fit compared to the original reference distributions. The improvement is particularly noticeable in the tails, where the residuals from the $\mathbb{H}$-extreme models exhibit a closer alignment with the theoretical reference lines. This behavior suggests that the proposed model more accurately captures the extreme-value behavior and tail dependence present in the data. In particular, visually the results from LSE seem a bit better than MLE, which is explained in the scenarios of small and moderate sample sizes. In those cases, MLE can suffer from high variance in the tail estimates because it relies heavily on the likelihood contributions of few extreme observations. LSE, which spreads the influence across all observations, can occasionally yield smoother (though less theoretically efficient, as demonstrated in Table \ref{tab:tab_estimation_results}) estimates in these regions. This does not affect our general conclusion that the $\mathbb{H}$-extreme distribution does provide a superior fit in the cases studied.}

\section{Conclusions}

We studied a general class of probability distributions with positive support, which generalizes several well-known and established models in the literature. In this study, we described various theoretical properties of this class. Additionally, to define the model's distribution function, we introduced the concept of the incomplete extreme value { $\mathbb{H}$-function}.


Two estimation methods were proposed for the new model. These methods were implemented computationally, and their results were validated by modeling three real-world datasets. These datasets had previously been modeled in existing literature using particular cases of the new model—specifically, the Fréchet, Weibull, and Gamma distributions. Our study demonstrates that this new class of distributions, despite having more parameters than the basic models, provides a better fit to the data. This conclusion is supported by information criteria based on penalized likelihood (AIC, BIC, and EDC).


Despite its superior fit, the new model presents computational challenges. These arise from the need to estimate six parameters, and the potential non-identifiability of the model can sometimes hinder or limit the performance of numerical optimization algorithms. Therefore, using good initial guesses and carefully controlling the search space are essential for successful optimization.

Potential future work could involve enhancing the optimization methods employed in this context. For instance, one approach might be to implement global search metaheuristics that build upon local search heuristics while also exploring variations in the search regions.

	\paragraph*{Acknowledgements}
The authors thank the support provided by the University of Brasilia (UnB). The research was supported in part by CNPq, CAPES, and FAPDF grants from the Brazilian government.
The authors are grateful to  J. Monteiro and M. Oliveira for helpful discussions.
 
	\paragraph*{Disclosure statement}
	There are no conflicts of interest to disclose.

\bibliographystyle{abbrv}
\bibliography{sample}


\end{document}